%% file: arxiv_dissertation.tex
\documentclass[final]{ua-thesis}
\pdfoutput=1

\usepackage{pdfpages}
\usepackage{amsmath, amssymb}
\usepackage{url}
\usepackage{enumerate}
\usepackage[all]{xy}

\newcommand{\B}{\mathcal B}

\newcommand{\EE}{\mathbb E}
\newcommand{\F}{\mathcal F}

\newcommand{\N}{\mathbb N}
\newcommand{\PP}{\mathbb P}
\renewcommand{\O}{\mathcal O}

\newcommand{\R}{\mathbb R}

\newcommand{\V}{\mathcal V}

\newcommand{\Z}{\mathbb Z}

\renewcommand{\span}{\operatorname{span}}

\newcommand{\SPD}{\operatorname{SPD}}
\newcommand{\Sym}{\operatorname{Sym}}

\newcommand{\Briem}{B^{\mathrm{R}}}
\newcommand{\Beuc}{B^{\mathrm{E}}}

\newcommand{\E}{\mathrm{e}}				
\newcommand{\D}{\mathrm{d}}				
\newcommand{\sD}{\, \mathrm{d}}		
\newcommand{\oo}{\infty}
\newcommand{\Leb}{\operatorname{Leb}}

\newcommand{\var}{\operatorname{Var}}

\numberwithin{equation}{chapter}
\theoremstyle{definition}
\newtheorem{env_thm}{Theorem}[chapter]
\newtheorem{env_def}[env_thm]{Definition}
\newtheorem{env_lem}[env_thm]{Lemma}
\newtheorem{env_cor}[env_thm]{Corollary}
\newtheorem{env_pro}[env_thm]{Proposition}
\newtheorem{env_con}[env_thm]{Conjecture}
\newtheorem{env_cla}[env_thm]{Claim}

\newtheorem{env_ass}[env_thm]{Assumption}

\director{Jan Wehr}
\author{Tom LaGatta}
\title{Geodesics of Random Riemannian Metrics}		
\date{2010}

\begin{document}

\maketitle

\chapter*{Acknowledgments}


	There is a long list of people I would like to thank who have helped get me through six years of graduate school.  First and foremost, I thank my advisor Janek Wehr, whose passion and taste are what drew me to probability in the first place, and who has guided and supported me more than anyone else.  I also thank the other members of my committee, Dave Glickenstein, Tom Kennedy and Joe Watkins, as well as Daniel Ueltschi, who have all been so generous with their advice, mathematical or otherwise.  I thank Chuck Newman, who will be my supervisor at the Courant Institute next year, whose work has been a foundation for my own, and who has been so very generous with his time and help this past year.
	
	I thank the members of the administration of the mathematics department at Arizona:  Doug Ulmer, Ken McLaughlin and Tom Kennedy, the successive heads of the graduate program; Nick Ercolani and Bill McCallum, the successive heads of the department; and Michael Tabor, the current PI for the VIGRE grant.  The math department would not be one tenth what it is without the excellent leadership of these men.  I also thank the staff at the department, most notably Sandy Sutton, the departmental secretary for the graduate program, who single-handedly manages to keep at bay the tide of paperwork we graduate students have to face.  Thanks also go out to the rest of the staff, including Sylvia Anderson, Jerrie Bieberstein, Bob Borys, Tina Deemer, David Gonzalez, Annette Horn, Denise Ingram, Anne Keyl, Christa King, Karl Newell, Alex Perlis, Karen Schaffner, Andrew Tubbiolo, Brooke Zang and every one else who has kept things running smoothly behind the scenes so that we grad students can focus on what we need to.  I also thank the administration at the upper echelons of the College of Science, whose commitment to graduate student excellence is unshakable in these times of budget cuts and transformations:  Joaquin Ruiz, Elliott Cheu, and Gail Burd, and their outstanding staff members including Danielle Shirar and Bernadette Thomas.
	
	There are many great teachers at Arizona, and I have been lucky to have worked with many of them.  I thank Bruce Bayly, Vita Borovyk, Sunhi Choi, Scott Clark, Carl Devito, Lennie Friedlander, Donna Krawczyk, Doug Pickrell, Natalie Sandler, Deidre Smith, Laurie Varecka, Mariamma Varghese, Bill V\'elez, Joe Watkins and Janek Wehr, who have all helped shape me into a better teacher.  
		
	My research has been supported for many semesters under the NSF VIGRE grant at the University of Arizona (DMS-06-02173), and for one semester under NSF grant DMS-06-23941.  I thank Krzysztof Gaw\c{e}dzki and ENS-Lyon for being outstanding hosts while I was in Lyon for a semester.
	
	I am thankful for the Epic Caf\'e, Caff\`e Luce and Grill in Tucson, where I spent countless hours in graduate school working, enjoying life, then working some more.  I thank my many, many friends in Tucson, San Francisco, New York City, Austin, Vancouver, Lyon and other cities across the globe, who continue to remind me how beautiful and diverse life can be, and how much opportunity it has to offer.  In this list I include my lovely officemates Jordan Allen-Flowers, Josh Chesler, Yaron Hadad, Jeffrey Hyman, Selin Kalaycioglu, Anya Petersen, Mandi Schaeffer Fry, Brad Weir, and Mei Yin, as well as my best friend Ben Polletta and my girlfriend Emily Chambliss.
		
	Finally and most importantly, I thank my parents, Barbara and Tom LaGatta.  I have reached this point only by standing on their shoulders.


\chapter*{Dedication}



{\centerline{To Tucson.}}


\tableofcontents

\begin{abstract}

We introduce Riemannian First-Passage Percolation (Riemannian FPP) as a new model of random differential geometry, by considering a random, smooth Riemannian metric on $\R^d$.  We are motivated in our study by the random geometry of first-passage percolation (FPP), a lattice model which was developed to model fluid flow through porous media.  By adapting techniques from standard FPP, we prove a shape theorem for our model, which says that large balls under this metric converge to a deterministic shape under rescaling.  As a consequence, we show that smooth random Riemannian metrics are geodesically complete with probability one.

In differential geometry, geodesics are curves which locally minimize length.  They need not do so globally:  consider great circles on a sphere.  For lattice models of FPP, there are many open questions related to minimizing geodesics; similarly, it is interesting from a geometric perspective when geodesics are globally minimizing.  In the present study, we show that for any fixed starting direction $v$, the geodesic starting from the origin in the direction $v$ is not minimizing with probability one.  This is a new result which uses the infinitesimal structure of the continuum, and for which there is no equivalent in discrete lattice models of FPP.  

\end{abstract}

\chapter{Introduction}

		Standard first-passage percolation (Standard FPP) is a model of random geometry on the discrete lattice, famously introduced by Hammersley and Welsh \cite{hammersley1965fpp} in order to model fluid flow through porous media.  The model is simple to define: take the lattice $\Z^d$ and associate to each bond (edge) a random number, called the passage time. This induces a random metric on $\Z^d$, where the distance between two points is the infimum of passage times over all paths which connect the two points.  Our model, Riemannian first-passage percolation, is a continuum analogue of Standard FPP.  Instead of a random discrete metric on the lattice $\Z^d$, we consider a random Riemannian metric in the continuum $\R^d$, which again gives rise to a random distance function. Both the lattice and continuum models have a similar global geometric structure, but Riemannian geometry provides a rich local structure to our model.
		
		Our consideration of a random Riemannian metric is a novel approach not found in the differential geometry literature.  For large-scale properties which do not depend on the local structure of the metric, we are able to directly adapt techniques from Standard FPP to our model.  To do this, we discretize the plane into unit cubes, and consider a dependent FPP model on the lattice formed by their centers.  We exploit this strategy in proving a shape theorem in our article \cite{lagatta2009shape}, included in Appendix A:  large balls under this metric converge to a deterministic shape under rescaling.  We also show that the metric is almost surely geodesically complete.

		In this dissertation, we sketch a proof of a new result which exploits the infinitesimal structure of our model, and which is not available for lattice models.  We assume now that our random metric has a rotationally invariant distribution, so that the limiting shape is a Euclidean ball.  As in differential geometry \cite{lee1997rmi}, we define geodesics to be curves which locally minimize arc length under our Riemannian metric.  Geodesics need not globally minimize length:  on the sphere, for example, geodesics are great circles, which do not minimize length past antipodal points.  Geodesics are defined by a local condition:  given a point and a direction, we define a geodesic to be the solution to a certain ordinary differential equation.  The completeness of the metric guarantees that geodesics can be extended for all time.  Geodesics can also be defined by a global condition as the curves which minimize the distance between two points, though in that case they need not be unique (consider the geodesics which connect antipodal points on a sphere).
		
		Our main result is that globally length-minimizing geodesics are rare, and the following event holds with probability one:  starting at the origin, the set $\V \subseteq S^{d-1}$ of directions which results in minimizing geodesics has Lebesgue measure zero on the sphere.  This measure-zero property is not a technicality:  the set $\V$ is non-empty and we furthermore conjecture that it is uncountable.  We believe our proof of the main result is correct, though we still have some technical details to finish and plan to submit it for publication soon.

		The proof of the main result is detailed and is split into three separate sections in Chapter \ref{presentstudy}.  We consider the geodesic $\gamma$ starting at the origin in a fixed direction $v \in S^{d-1}$.  By adapting our techniques from \cite{lagatta2009shape}, we show in Section \ref{sect_frontiertimes} that there is a sequence of ``frontier times'' $t_k$ along the geodesic at which the metric is ``well-behaved'' in a neighborhood $B_k$ of $\gamma(t_k)$.  As the geodesic exists the Euclidean ball $\Beuc(0, |\gamma(t_k)|)$, being well-behaved means that its exit velocity satisfies a cone condition, and the geodesics are bounded uniformly away from being tangential to the ball.  The metric is well-behaved in the sense that we have a uniform bound on the $C^{2+\alpha}$-norm of $g$ in $B_k$, as well as a lower bound for the minimum eigenvalue of $g$ in $B_k$.  
				
		In Section \ref{sect_unifprob}, we show that there is a uniform lower bound $p > 0$ of the probability that a certain event $U_k$ occurs at the frontier time $t_k$.  In the proof, we change perspectives from the FPP context and focus on probability measures on Banach spaces, by means of a Strong Markov Property.  The Arzel\`a-Ascoli theorem \cite{folland1999real} implies that the set $\Gamma$ of well-behaved metrics on $B_k$ is compact in the space of $C^2$-metrics.  We show that conditional probabilities in this context vary continuously on the conditioning, so by minimizing over the compact set $\Gamma$ we have a positive lower bound on the probability of $U_k$. In order to carry out this uniform probability argument, the author developed the concept of continuous disintegrations, which evolved into the separate publication \cite{lagatta2010continuous}, included in Appendix B.  
			
		It is the presence of positive curvature which destabilizes minimizing geodesics \cite{lamburt2003grc}.  In Section \ref{sect_bump}, we exploit this observation, and we describe a way to extend the metric at $\gamma_v(t_k)$ to a bump metric ahead, and argue that geodesics cannot be minimizing after spending enough time on the bump (like the top two-thirds of a sphere).  This property is perturbed under small perturbations of the bump metric, so the event $U_k$ is that the random metric $g$ is sufficiently close to the bump metric in the region ahead.  
		
		Finally, we put the pieces together to prove the main result.  We fix a direction $v \in S^{d-1}$, and estimate the probability that $v \in \V$ (i.e. that the geodesic $\gamma$ which it generates is minimizing).  If there is not a sequence of frontier times $t_k$ as described above, the geodesic is not minimizing; supposing there is such a sequence, at each $t_k$ there is a uniform probability $p$ that the geodesic runs over a bump and stops being minimizing.  Consequently, with probability one, $v \notin \V$.\newline
		
	For the remainder of this Introduction, we give a review of the literature for Standard FPP (Section \ref{standardfpp}) and related models (Section \ref{models}).  
	In Section \ref{includedpapers}, we described the articles which we include in Appendices A and B.  
	
	[n.b.: This version of the dissertation does not include the two referenced papers \cite{lagatta2009shape} and \cite{lagatta2010continuous}.]

\section{Standard First-Passage Percolation} \label{standardfpp}

	We formally introduce the model of Standard FPP.  Consider the $d$-dimensional lattice $\Z^d$ with $d \ge 2$.  Let $\{t_b\}$ be a family of independent, identically distributed, non-negative random variables, indexed by bonds (nearest-neighbor edges) $b$ of the lattice.  For any $z,z' \in \Z^d$, define the passage time from $z$ to $z'$ by
		$$\tau(z,z') = \inf_{\gamma} \sum_{b \in \gamma} t_b,$$
	where the infimum is taken over all lattice paths $\gamma$ connecting $z$ to $z'$.  This $\tau$ is a random distance function on $\Z^d$.  For a very good introduction to Standard FPP, see Howard \cite{howard2004mfp} or the more recent Blair-Stahn \cite{blair2007first}.  
	
\subsection{Time Constant}
	
	The first object of study is the passage time $a_n := \tau(0, n\E_1)$ between the origin and the point $n\E_1 = (n,0,\dots,0)$.  One wishes to study the asymptotic behavior of this quantity as $n \to \infty$.  In \cite{kingman1968ets}, Kingman formulated his famous subadditive ergodic theorem in order to prove the basic result of FPP:  
	
	\begin{env_thm}
	If the passage times have finite mean, there exists a non-random constant $\mu_{\E_1}$, such that 
		$$\lim_{n\to\oo} \tfrac{1}{n} \tau(0,n\E_1) = \mu_{\E_1}$$
	almost surely and in $L^1$.  
	\end{env_thm}
	
	By symmetry, the same value is the limit for the passage times $\frac 1 n \tau(0,n\E_i)$ in any of the coordinate directions $\E_i$.  More generally, for each direction $v \in S^{d-1}$, there exists a non-random constant $\mu_v$, such that
		\begin{equation} \label{kingmandiscrete}
			\lim_{n\to \oo} \tfrac{1}{n}\tau(0,\widetilde{nv}) = \mu_v \end{equation}
	almost surely and in $L^1$, where $\widetilde{nv} \in \Z^d$ is the nearest lattice point to $nv$.  The constants $\mu_v$ vary continuously with $v$ (cf. \cite[Proposition 1.3]{lagatta2009shape}).  Kesten \cite{kesten1180arp} has shown that the constant $\mu_v$ is non-zero provided that the probability that $t_b = 0$ is less than the critical percolation probability for $\Z^d$; see \cite{howard2004mfp} for more details.  While Kingman's theorem asserts the existence of the time constants $\mu_v$, they have not been computed explicitly for any non-trivial distribution of passage times for Standard FPP.

\subsection{Shape Theorem}

	Henceforth, we assume that $\mu_v > 0$ for all $v \in S^{d-1}$, and that the passage time distribution satisfies the simple moment condition 
		\begin{equation} \label{shapecond}
			\EE \min\{t_1, \dots, t_{2d}\}^{2d} < \oo, \end{equation}
	for $2d$ independent copies $t_1, \cdots, t_{2d}$ of $t_b$.  Where \eqref{kingmandiscrete} is a law of large numbers-type statement for each fixed direction $v$, the shape theorem of Cox and Durrett \cite{cox1981slt} is a stronger result which holds for all directions simultaneously.  We extend the function $\mu_v$ to a norm on $\R^d$ by defining $\mu(x) := \mu_{x/|x|} |x|$.  Consider the unit ball in this norm,
		$$A = \{ x : \mu(x) \le 1 \} = \{x : |x| \le \mu^{-1}_{x/|x|} \}.$$
	This non-random set depends only on the distribution of the passage times $t_b$, and is convex, compact, and invariant under the symmetries of the lattice $\Z^d$.
	
	Consider
		$$\tilde B_t = \{z \in \Z^d : \tau(0,z) \le t\},$$
	the random ball of radius $t$ in $\Z^d$.  This is a lattice object, so we ``inflate'' it to get a continuum one:  for $z \in \Z^d$, let $C_z = [z-1/2,z+1/2)^d$ be the unit cube centered at $z$ in $\R^d$, and let 
		$$B_t = \bigcup_{z \in \tilde B_t} C_z.$$
	We define the rescaling $\tfrac{1}{t} B_t$ as the set of all points $x \in \R^d$ such that $tx \in B_t$.  The shape theorem says essentially that $\tfrac{1}{t} B_t \to A$:
	
	\begin{env_thm}
	For all $\epsilon > 0$, with probability one, there exists a time $T$ such that if $t \ge T$, then
		\begin{equation} \label{shapeformula}
			(1-\epsilon)A \subseteq \tfrac{1}{t} B_t \subseteq (1+\epsilon)A. \end{equation}
	\end{env_thm}
	
	While the existence of a limiting shape $A$ is guaranteed by this theorem, it is in practice and in theory very difficult to obtain much information on the precise shape $A$.  
	For Standard FPP, there are no known passage-time distributions which yield a rotationally-invariant limiting shape---lattice effects always seem to persist \cite{howard2004mfp}.   Durrett and Liggett \cite{durrett1981shape} have shown that if the distribution of $t_b$ has a positive atom with sufficiently high probability, then there are ``facets'' in the limiting shape:  where it meets the diagonal directions, $\partial B$ is made up of flat pieces.  
	
	One expects that the facets of Durrett and Liggett are pathological, and that the boundary $\partial A$ should typically satisfy some smoothness properties.  Newman and Piza \cite{newman1995dsf} define a direction of curvature $v \in S^{d-1}$ if $\partial A$ is locally spherical near $x := v/\mu_v \in \partial A$.  Precisely, this means that there exists a Euclidean ball $D$ depending on the direction $v$ which contains $A$ and is tangent to $A$ at $x$:
		$$A \subseteq D \mathrm{~and~} x \in \partial D.$$

	There is a simple proof that directions of curvature exist in Standard FPP:  let $r$ be the minimal radius such that the Euclidean ball $D = B(0,r)$ centered at the origin contains $A$, then the directions where $D$ meets $A$ are directions of curvature \cite{howard2004mfp}.  This is a very weak existence result, however, and there may only be finitely many directions of curvature.  Moreover, no specific direction has been verified to be a direction of curvature for any distribution of passage times, including the axial directions \cite{howard2004mfp}.  
	
	Newman \cite{newman1995svf} says that $A$ is uniformly curved if every direction is a direction of curvature, and moreover that the radii of the balls $D$ is uniformly bounded away from infinity.  In the $d=3$ case where $\partial A$ is a topological $2$-sphere, this assumption is that the Gaussian curvature is uniformly bounded away from $0$ on the boundary surface $\partial A$.  Again, this has not been verified for any particular distribution of passage times in Standard FPP.  However, in the rotationally invariant models of Euclidean FPP of Newman and Howard \cite{howard1997euclidean} and Riemannian FPP of LaGatta and Wehr \cite{lagatta2009shape}, the limiting shape is a Euclidean ball.  Consequently, all directions are directions of curvature, and the limiting shape is uniformly curved.

\subsection{Shape Fluctuations and $\chi$} \label{shapefluct}
	
	As is to be expected from a law of large numbers, the upper bound $\epsilon t$ in \eqref{shapeformula} on the fluctuations of $B_t$ from $tA$ is far from optimal.  Using an exponential moment condition, Kesten \cite{kesten1993scf} was able to improve \eqref{shapeformula} to
	\begin{env_thm}
		$$(t-ct^{\kappa} \log t) A \subseteq B_t \subseteq (t+ct^{1/2} \log t) A$$
	for some non-random constant $c > 0$ and $\kappa < 1$.  	
	\end{env_thm}
	
	In his proof, he also showed that given a second-moment condition, the variance of the passage time $a_n = \tau(0,n\E_1)$ was at worst linear:
		\begin{equation} \label{linvar}
			\var a_n \le C n, \end{equation}	
	for some non-random constant $C > 0$.  
	
	Let us define the longitudinal fluctuation exponent $\chi$ as the minimum number $k$ such that, with probability one, there exists a time $T$ such that if $t \ge T$, then
		\begin{equation} \label{chidef}
			(t-t^k) A \subseteq B_t \subseteq (t+t^k) A, \end{equation}
	Alexander \cite{alexander1993note, alexander1997approximation} was able to improve Kesten's estimate and proved \eqref{chidef} with the value $k = 1/2$.  Consequently, in terms of the fluctuation exponent this is the upper bound $\chi \le 1/2$.
	
	The above discussion may suggest that $a_n$ behaves diffusively.  In fact, this is far from the case.  First-passage percolation is conjectured to lie in the Kardar-Parisi-Zhang (KPZ) universality class of growth processes \cite{kardar1986dynamic,krug1988universality, krug1991kinetic}.  
	In two-dimensions, the optimal value of the fluctuation exponent $\chi$ should be $1/3$, and it is believed that the variance $\var a_n$ is of order $n^{2/3}$.  We will explore this connection in more detail in the next section. 

	Benjamini, Kalai and Schramm \cite{benjamini2003first} used the concentration inequalities of Talagrand \cite{talagrand1994russo} in order to show that $a_n$ has sublinear distance variance
		$$\var a_n \le C n / \log n$$
	for a Bernoulli distribution of passage times, and Bena\"im and Rossignol \cite{benaim2006modified, benaim2008exponential} were able to extend this to a much wider class of passage time distributions, including exponential passage times.  This may seem like a trivial improvement of \eqref{linvar}, but in fact is quite significant.  Chatterjee \cite{chatterjee2008chaos} has some remarkable applications of sublinear distance variance (or ``superconcentration'' in his terminology) which have not yet been applied successfully to first-passage percolation; he has, however, done this for spin-glass models  \cite{chatterjee2009disorder}, and we will discuss this more in Section \ref{spinglasses}.
	

\subsection{Transversal Fluctuations and $\xi$} \label{KPZ}

	If the passage time distribution does not have any atoms, then with probability one, for all $n$ there exists a unique path $\gamma_n$ which realizes the minimum passage time $a_n = \tau(0,n\E_1)$ from the origin to $n\E_1$ \cite{howard2004mfp}.  Let $d_n$ be the maximal distance that $\gamma_n$ deviates from the straight line path from $0$ to $n\E_1$.  Formally,
		$$d_n = \sup \left\{ \inf_{0 \le j \le n} |\gamma_n(i) - j\E_1| : 0 \le i \le |\gamma_n| \right\},$$
	where $|\gamma_n|$ is the number of points in the path $\gamma_n$, and $|\gamma_n(0) - j\E_1|$ is the Euclidean distance between the points $\gamma_n(i)$ and $j\E_1$ in $\Z^d \subseteq \R^d$.

	We define the transversal fluctuation exponent $\xi$ as the minimum number such that with probability one,
		$$d_n = O(n^\xi).$$
	We can similarly define the exponent $\xi(v)$ in the direction $v \in S^{d-1}$, though it is believed that the quantities $\xi(v)$ are invariant under direction.  As mentioned in the previous section, it is conjectured that first-passage percolation lies in the KPZ universality class of growth processes \cite{krug1991kinetic}, and the fluctuation exponents $\chi$ and $\xi$ satisfy the KPZ equation
		\begin{equation} \label{KPZeqn}
			\chi = 2\xi - 1. \end{equation}
	It is further conjectured in dimension $d=2$ that $\chi = 1/3$ and $\xi = 2/3$ \cite{krug1991kinetic}.	
	
	In terms of rigorous results, Newman and Piza \cite{newman1995dsf} have partially proved an inequality of the form \eqref{KPZeqn}, but only as an inequality, and for transversal fluctuations in directions of curvature: 
	
	\begin{env_thm}
	If $v$ is a direction of curvature,
		\begin{equation} \label{NPupper}
			\chi \ge 2\xi(v) - 1. \end{equation}
	\end{env_thm}
	Along with the Kesten-Alexander upper bound $\chi \le 1/2$, this implies the upper bound $\xi(v) \le 3/4$ on transversal fluctuations in directions of curvature.  Using techniques based on Wehr-Aizenman \cite{wehr1990fluctuations}, Newman and Piza \cite{newman1995dsf} are also able to prove the lower bound
		\begin{equation} \label{NPlower}
			2\chi'(v) \ge 1 - (d-1) \xi(v) \end{equation}
	for an exponent $\chi'(v)$ related to $\chi$.  They conjecture that $\chi'(v)$ is independent of direction and is in fact equal to $\chi$. 
	
	If the KPZ equation \eqref{KPZeqn} holds, then the trivial bound $\chi \ge 0$ implies that $\xi \ge 1/2$.  This is non-trivial, since the value $\xi = 1/2$ corresponds to the process $d_n$ behaving diffusively (e.g. like a simple random walk).  However, it is believed that $d_n$ behaves super-diffusively, and $\xi > 1/2$.  Under some weak conditions, Licea, Newman and Piza \cite{licea1996superdiffusivity} rigorously prove  the lower bound $\xi \ge 1/(d+1)$ for all dimensions $d$, as well as $\xi'(d) \ge 1/2$ for a related exponent $\xi'(d)$ depending on dimension $d$, and $\xi'(2) \ge 3/5 > 1/2$.  It is conjectured that $\xi' = \xi$, but it is still an open question to prove rigorously that $\xi > 1/2$ for any model of FPP.

\subsection{Geodesics and Disordered Ferromagnets} \label{gadf}

	As evidenced from the above section, minimizing paths are of critical important to the study of first-passage percolation:  the fluctuations of minimizing paths are related to the fluctuations of the limiting shape via the KPZ equation \eqref{KPZeqn}.  Assume that the passage time distributions are continuous.  This implies that finite minimizing paths exist; however, the existence of infinite minimizing paths is a subtler question.
	
	In the first-passage percolation literature, minimizing paths are denoted by the term ``geodesic.''  This is very different from the standard meaning of the word in differential geometry.  As we will discuss in more detail in Chapter \ref{presentstudy}, geodesics are curves which locally minimize length, but not necessarily globally.  On the sphere, for example, geodesics are great circles, which do not globally minimize distance past antipodal points.  In this section, we only use the term geodesic to refer to minimizing lattice paths, since there is not an infinitesimal notion of geodesic for lattice models.  However, when we discuss the model of Riemannian FPP in Chapter \ref{presentstudy}, we will distinguish between ``geodesics'' which only locally minimize length, and ``minimizing geodesics'' which do so globally.
	
	We say that a path $\gamma : \N \to \Z^d$,
		$$\gamma = (\gamma^0, \gamma^1, \gamma^2, \dots ),$$
	is a one-sided geodesic if for every pair of points $x,y \in \gamma$, the passage time $\tau(x,y)$ is realized as the passage time along $\gamma$.  Provided that the passage time distribution has no atoms, between any two points there exists a unique minimizing path \cite{howard2004mfp}.  It is easy to extend this using a spanning-tree argument to show that with probability one, for every point $z \in \Z^d$, there exists a one-sided geodesic starting at $z$.  Fix $z$, and for every point $z'$ let $\gamma_{z'}$ denote the unique minimizing path connecting $z$ and $z'$.  Let
		$$\mathcal T(z) = \bigcup_{z'} \gamma_{z'}$$
	be the union of the edges of these minimizing paths.  Clearly, $\mathcal T(z)$ is a spanning tree of $\Z^d$ hence contains an infinite path starting at $z$.
	
	The above demonstrates that there is at least one geodesic at each point, though Newman \cite{newman1995svf} conjectures that there should be infinitely many.  For $w \in S^{d-1}$, we say that $w$ is an asymptotic direction for a geodesic if the limit 
		$$\lim_{n\to\oo} \frac{\gamma^n}{|\gamma^n|}$$
	exists and equals $w$.  Under the assumption of uniform curvature on the limiting shape, Newman shows that, with probability one, every one-sided geodesic at the origin has an asymptotic direction.  Furthermore, every direction $w \in S^{d-1}$ is realized as the asymptotic direction for at least one geodesic, which implies that there are infinitely many geodesics at the origin.  While the uniform curvature assumption is not satisfied for any known distribution of passage times, these arguments have been successfully applied to other models, which we discuss more in Section \ref{models}.
	
	For Standard FPP, H\"aggstr\"om and Pemantle \cite{häggström1998first} are able to show that if $d=2$ and the passage times have an exponential distribution, then with positive probability any particular site (e.g. the origin) has at least two distinct one-sided geodesics.  Their argument involves a connection to Richardson's growth model \cite{richardson1973random}.  Hoffman \cite{hoffman2008geodesics} extends their ideas to show that the number of one-sided geodesics at the origin is at least $4$ with positive probability.  The number $4$ comes from the minimum number of sides of the limiting shape $A$ in $d=2$.  If $A$ is not polygonal (for example a ball), then we say it has infinitely-many sides, and accordingly the number $4$ is improved to $\oo$.  Unfortunately, like the uniform curvature assumption, that $A$ is not polygonal has not been rigorously shown for any passage time distribution.
	
	We extend the definition of geodesic to two-sided paths $\gamma : \Z \to \Z^d$,
		$$\gamma = (\dots, \gamma^{-1}, \gamma^0, \gamma^1, \dots ),$$
	and say that $\gamma$ is a two-sided geodesic if for every pair of points $x,y \in \gamma$, the passage time $\tau(x,y)$ is realized as the passage time along $\gamma$.  The existence of two-sided geodesics is an important open question, and it is believed that the answer is different for dimensions $d=2$ and $d > 2$.  In either case, Wehr \cite{wehr1997number} has shown that if two-sided geodesics do exist, then there are infinitely many of them with probability one.  
	
	The existence of two-sided geodesics has consequences for statistical physics.  In $d=2$, Standard FPP is essentially the dual model for a disordered ferromagnet, a simplification of the Edwards-Anderson spin glass \cite{edwards1975theory} where the nearest-neighbor couplings are non-negative random variables \cite{newman1997tds}.  Consequently, the almost-sure existence of two-sided geodesics is equivalent to the almost-sure existence of non-trivial ground states in this model.  For physical reasons, it is conjectured that these do not exist.  We will discuss spin-glass models more in Section \ref{spinglasses}; for more details on the connection between these two models, see Newman \cite{newman1997tds}.
	
	Under the assumption of uniform curvature on the limiting shape $A$, Newman \cite{newman1995svf} shows that the only two-sided geodesics which can exist are those with antipodal asymptotic directions.  i.e., with probability one, for each two-sided geodesic $\gamma$ there exists $w \in S^{d-1}$ such that
		$$\lim_{n\to \pm\oo} \frac{\gamma^n}{|\gamma^n|} = \pm w.$$
	In the $d=2$ case, Licea and Newman \cite{licea1996gtd} show that for each deterministic $w$ that, with probability one, there does not exist a two-sided geodesic with asymptotic directions $w$ and $-w$.  Their argument is to fix $w \in S^{d-1}$, and consider the minimizing paths $\gamma_n$ from $-nw$ to $nw$ (properly adjusted so the points $\pm nw$ lie on the lattice).  If a two-sided geodesic $\gamma$ with asymptotic directions $\pm w$ were to exist, then $\gamma_n$ should converge to it; if $d_n$ equals the minimal distance from the origin to $\gamma_n$, then $d_n$ is of order $1$.  However, by Section \ref{KPZ}, $d_n$ should scale like $n^\xi$, so a positive lower bound on $\xi$ gives a contradiction.
	
	One of the few solid non-existence results is due to Wehr and Woo \cite{wehr1998absence}, who have shown that for FPP restricted to the half-lattice in $d=2$, there exist no two-sided geodesics with probability one.

	There is a heuristic scaling argument that suggests superdiffusivity of transversal fluctuations (i.e. $\xi > 1/2$) implies non-existence of two-sided geodesics in $2$-dimensional Standard FPP \cite{newmanemailapr22}.  Consider a circle of large radius $R$ centered at the origin, and break it into arcs of length $O(R^\xi)$, so that there are $O(R^{1-\xi})$ such arcs $\alpha_i$.  If two-sided geodesics exist with probability one, then one meets the origin with some probability $p > 0$ not depending on $R$.  By coalescence arguments \cite{licea1996gtd, howard1997euclidean}, for all points $x \in \alpha_i$ and $y \in -\alpha_i$ on antipodal arcs, the minimizing paths from $x$ to $y$ should all coalesce with high probability.  Thus there is essentially only one geodesic passing between antipodal arcs, whose transversal fluctuations are of order $O(R^{\xi})$ by definition of $\xi$; consequently, the probability that it passes through the origin is of order $O(R^{-\xi})$.  Roughly speaking, for different arcs, these events are almost independent, so by considering them as Bernoulli trials, the probability of at least one occurring is $O(R^{1-\xi} \cdot R^{-\xi}) = O(R^{1-2\xi})$.  If $\xi > 1/2$, this probability goes to zero; in particular, for large $R$ it is less than $p$, a contradiction.






\section{Other Models Related to First-Passage Percolation} \label{models}

	Boivin \cite{boivin1990first} introduced a model of stationary, ergodic FPP, where one considers $d$ different passage time distributions, one for each direction $\E_1, \dots, \E_d$, and where the assumption of independence of passage times is relaxed to ergodicity.  He shows that a shape theorem is satisfied, and that the limiting shape $A$ in this setting need only be compact, convex and satisfy the antipodal symmetry $A = -A$.  H\"aggstr\"om and Meester \cite{haggstrom1995ass} show that these conditions on $A$ are sufficient for there to exist a stationary, ergodic passage time distribution with limiting shape $A$.  Consequently, there exists some distribution of stationary, ergodic passage times for which $A$ is a Euclidean ball.

	Chatterjee-Dey \cite{chatterjee2009central} consider a model of first-passage percolation restricted to growing cylinders of the form 
		$$[0,n] \times [-h_n, h_n] \subseteq \Z^2$$
	where $h_n = o(n^{1/3})$ (as well as the generalization to higher dimensions).  In this context, they show that the first passage times $a_n = \tau(0,n\E_1)$ satisfy a Gaussian central limit theorem; in particular, the fluctuations of $a_n$ are order of $\sqrt{n}$.  This is qualitatively different than the expected behavior in Standard FPP, where it is believed that $a_n$ has fluctuations of order $n^\chi$, and does not satisfy a Gaussian central limit theorem.  Thus this is some rigorous evidence suggesting the lower bound $\chi \ge 1/3$ for Standard FPP.

\subsection{Euclidean First-Passage Percolation}

	In Standard FPP, lattice effects always seem to persist at the macroscopic scale:  the limiting shape $A$ is not rotationally-invariant for any known passage time distribution.  As we saw in the last section, this is a major obstruction, as many results for Standard FPP on fluctuation exponents or the existence of geodesics hold only under strict curvature assumptions on the limiting shape.  In order to circumvent this rigidity, Vahidi-Asl and Wierman \cite{vahidi1990first,vahidi1992shape} consider a model of FPP on the Voronoi graph generated by the points of a random Poisson point process $Q$ in $\R^d$.  They are able to prove a shape theorem and, since $Q$ has a rotationally-invariant distribution, the limiting shape is a Euclidean ball.  The shape theorem for this model is non-trivial for any distribution of passage times---including the case when passage times are non-random and constant---as the Voronoi graph of a random point process may be quite complicated.	
	
	Howard and Newman \cite{howard1997euclidean} consider another model based on a Poisson point process $Q$ in $\R^d$, which they term Euclidean FPP.  Rather than deal with the complicated spatial structure of the Voronoi graph, they work with the complete graph on $Q$, where every point is adjacent to each other.  They fix a parameter $\alpha > 1$, and define the passage time $t_b$ of the bond $b$ connecting the lattice points $q$ and $q'$ to be
		$$t_b = |q - q'|^\alpha,$$
	where $|\cdot|$ denotes the Euclidean norm in $\R^d$.  By the triangle inequality of the norm, if $0 \le \alpha \le 1$ then the passage time between $q$ and $q'$ is trivially minimized by taking the single edge between them.  However, when $\alpha > 1$, long jumps are discouraged and the model is non-trivial.
	
	Howard and Newman \cite{howard1997euclidean} prove a shape theorem with limiting shape a Euclidean ball.  Consequently, many of the results presented for Standard FPP in Section \ref{standardfpp} under restrictive hypotheses hold automatically for Euclidean FPP in for all $\alpha > 1$.  In \cite{howard2000geodesics}, Howard and Newman prove the shape fluctuation exponent bound $\chi \le 1/2$ using a moderate-deviations estimate similar to Kesten's \cite{kesten1993scf}, as well as the inequality \eqref{NPupper} so that $\xi \le 3/4$.  Howard \cite{howard2000lower} proves the inequality \eqref{NPlower} and the lower bound $\xi \ge 1/(d+1)$.

	
	For geodesics, Howard and Newman \cite{howard2000geodesics} show that with probability one, every one-sided geodesic has an asymptotic direction, and for every $w \in S^{d-1}$ there exists a one-sided geodesic with asymptotic direction $w$.  In the $\alpha \ge 2$, they have slightly stronger results \cite{howard1997euclidean}, which they believe should also hold for the $\alpha > 1$ case.  As with other models of FPP, the existence of two-sided geodesics is still open.  The strongest theorem on two-sided geodesics in this setting \cite{howard2001special} is as with Standard FPP:  with probability one, all two-sided geodesics (if they exist) must have antipodal asymptotic directions; and that for any deterministic $w \in S^{d-1}$, with probability one there do not exist any two-sided geodesics with asymptotic directions $\pm w$.  
	

\subsection{Last-Passage Percolation} \label{LPP}


	Consider the two-dimensional lattice restricted to the upper-right quadrant,
		$$\Z^2_+ = \left\{ (z^1, z^2) \in \Z^2 : z^1 \ge 0 \mathrm{~and~} z^2 \ge 0 \right\}.$$
	We assign i.i.d. non-negative passage times $t_b$ to each bond $b$, and define the last-passage times between $z$ and $z'$ in $\Z^2$ to be
		$$\tau(z,z') = \sup_{\gamma} \sum_{b \in \gamma} t_b,$$
	where the supremum is taken over  directed paths $\gamma$ moving to the up or to the right.  The last-passage time is superadditive,
		$$\tau(z,z') \ge \tau(z,w) + \tau(w,z'),$$
	rather than subadditive like in FPP, but a superadditive version of Kingman's ergodic theorem \cite{martin2000linear} can be applied to prove the existence of a time constant, and similarly a shape theorem \cite{martin2004limiting}.  
	
	Directed FPP models are similar to directed polymer growth in physics, where the role of passage times is replaced by random potentials; we explore this connection further in the next section.  The first-passage time between two points represents the minimal energy of a polymer, and the last-passage time the maximal energy.	
	
	Even though the models of FPP and LPP are qualitatively different, they are both believed to lie in the KPZ universality class \cite{hambly2007heavy} and consequently satisfy the KPZ equation \eqref{KPZeqn}.  Most impressively, many exact results (scaling laws, asymptotic distributions, have been found for LPP, mostly for exponentially or geometrically distributed passage times \cite{hambly2007heavy}.  In the case of geometric passage times, Johansson \cite{johansson2000shape} exploits a beautiful connection to random matrix theory by means of increasing subsequences of random permutations \cite{baik1999distribution}.  Consider the last-passage time $a_n = \tau(0,n\E_1)$.  Johansson explicitly computes the time constant $\mu$ for which $\tfrac{1}{n} a_n \to \mu$ and its variance $\sigma^2 n^{2/3}$, implying that the shape fluctuation exponent $\chi$ is exactly equal to $1/3$.  More incredibly, he shows that $a_n$ satisfies a central limit theorem, and that the normalized random variable
		$$ \frac{a_n - \mu n}{\sigma n^{1/3}}$$
	converges in distribution to the Tracy-Widom distribution \cite{tracy1994level}.
	

\subsection{Directed Polymers in a Random Environment} \label{dirpolys}

	We describe a model of directed polymer growth in the presence of random impurities, following the survey By Comets, Shiga and Yoshida \cite{comets2004probabilistic}.  This is a particular model of random walk in a random environment which in spirit has many similarities to first-passage percolation.  In particular, in a certain temperature regime, the shape and transversal fluctuations of polymers are believed to satisfy the KPZ equation \eqref{KPZeqn}.  The literature on random polymer models is vast; we recommend the books by Giacomin \cite{giacomin2007random} and den Hollander \cite[Chapter 12]{den2009random}.
	
	Fix $d \ge 1$, and consider the state space $\N \times \Z^d$.  A polymer is a randon path $\{(j,\omega_j)\}_{j=1}^n$ in this space, increasing deterministically in the time coordinate $j$.  In the absence of impurities, the distribution of the path $\omega_j$ will a simple random walk in $\Z^d$ starting from the origin.  The effects of the impurities are summarized by random variables $\eta(n,x)$ at each site of $\N \times \Z^d$.  These variables $\eta$ represent random potential energies, and polymers will tend toward sites where $\eta$ is positive.  We scale the energies by the non-negative parameter $\beta$.  As we will see, there is a phase transition in $\beta$, which depends on the dimension $d$.
	
	We present the model formally.  Let $(\Omega_0, \F_0, \PP)$ be a probability space, and let $\eta = \{\eta(n,x)\}$ be a family of i.i.d. random variables on $\Omega_0$ indexed by $\N \times \Z^d$.  The value $\eta(n,x)$ represents the potential energy at time $n \in \N$ at site $x \in \Z^d$.  Write $\EE$ for expectation with respect to $\PP$, and suppose that $\eta(n,x)$ has a finite moment-generating function:
		$$\EE[ \E^{\beta \eta(n,x)} ] < \oo$$
	for all $\beta \in \R$.
		
	Write $\Omega$ for the space of paths $\omega = \{\omega_j\}$ in $\Z^d$ starting at the origin.  For $n > 0$, define the random Hamiltonians $H_n : \Omega \to \R$ by
		$$H_n(\omega) = - \sum_{j=1}^n \eta(j,\omega_j).$$
	Let $\F$ be the cylinder $\sigma$-algebra on $\Omega$, and let $\nu$ be the probability measure on $(\Omega, \F)$ so that under $\nu$, $\omega \in \Omega$ is a simple random walk on $\Z^d$ starting from the origin.  
	
	The measure $\nu$ is a background measure on the space of paths $\Omega$, which we modify by the effect of the environment.  Let $\beta \ge 0$ be a non-negative parameter, and define the random Gibbs measure $\mu_n$ on $\Omega$ by
		$$\mu_n(\omega) = Z_n^{-1} \E^{-\beta H_n(\omega)} \sD \nu(\omega),$$
	where the normalizing constant $Z_n = \int_\Omega \E^{-\beta H_n(\omega)} \sD \nu(\omega)$ is called the partition function of $\mu_n$.  An important quantity is the free energy
		$$F_n = -\tfrac{1}{\beta} \log Z_n.$$
	The Gibbs measure $\mu_n$, the partition function $Z_n$ and the free energy $F_n$ are all random with respect to the probability measure $\PP$.  
	
	The parameter $\beta$ represents the inverse temperature (precisely, $\beta^{-1}$ equals the temperature multiplied by the Boltzmann constant).   When $\beta = 0$, the system is in the infinite-temperature regime, and polymers are exactly simple random walks.  When $\beta \approx 0$ (high temperature), if polymers behave similarly to simple random walk, we say the system is in the weak-disorder phase.  Conversely, when $\beta \to \oo$, the thermal fluctuations of polymers is suppressed, and they should behave like minimizing paths in first-passage percolation.  We call this the strong-disorder phase.  Things are qualitatively quite different in the two regimes:  as we have seen in the previous sections, minimizing paths in FPP are superdiffusive, whereas simple random walk is diffusive.  

	The phase transition depends on the dimension $d$:    if $d = 1$ or $2$, the system is in the strong disorder phase in the presence of any disorder (i.e. $\beta \ne 0$); in the case $d \ge 3$, there exists a non-trivial critical value $\beta_c > 0$ such that when $\beta < \beta_c$ the system is in the weak disorder phase, and when $\beta > \beta_c$ it is in the strong disorder phase \cite{comets2004probabilistic}.  Derrida \cite{derrida1990directed} has estimates on the critical value $\beta_c$ in terms of the dimension $d$.
	
	The model was introduced in the physics literature by Huse and Henley \cite{huse1985pinning} in order to model interface boundaries in the low-temperature regime of an Ising model with random impurities.  In the $d=1$, strong disorder phase (which corresponds to 2-dimensional FPP), they gave numerical evidence that the transversal fluctuations of polymers should scale like $n^\xi$ for $\xi = 2/3$.  Soon after, this value for $\xi$ was confirmed by Huse, Henley and Fisher \cite{huse1985huse} and Kardar and Nelson \cite{kardar1985commensurate} using heuristic, physical arguments.  
	
	In Section \ref{shapefluct}, we saw that the shape fluctuation exponent $\chi$ is related to the variance of the passage time $a_n = \tau(0,n\E_1)$ in Standard FPP.  Here, the role of passage times is replaced by energy, and $\chi$ is related to the fluctuations of the free energy $F_n$:
		$$\var F_n \sim n^{2\chi},$$
	where the variance is with respect to the probability measure $\PP$.  It is strongly believed \cite{krug1991kinetic} that this model also lies in the KPZ universality class, hence satisfies the KPZ equation
		$$\chi = 2\xi - 1$$
	in all dimensions $d$.  When $d = 1$, it is believed that $\chi = 1/3$, as with two-dimensional Standard FPP.
	
	Imbrie and Spencer \cite{imbrie1988diffusion} formulated the mathematical model described above; for the Bernoulli potential $\eta = \pm 1$, they rigorously showed the phase transition in $d \ge 3$ using an expansion in the small parameter $\beta$.  Bolthausen \cite{bolthausen1989note} reproved their result using a simple martingale method, which Song and Zhou \cite{song1996remark} extended for general environments $\eta$.  Adapting the uniform-curvature assumption for Standard FPP as described in Sections \ref{shapefluct} and \ref{KPZ}, Piza \cite{piza1997directed} proved many rigorous results on fluctuation exponents in this model.


	When $d \ge 3$ and $\beta$ is small (weak-disorder phase), one expects polymers to behave roughly like simple random walks.  Carmona and Hu \cite{carmona2002partition} proved a theorem on delocalization of polymers in a Gaussian random potential $\eta$ (later improved to general potentials by Comets, Shiga and Yoshida \cite{comets2003directed}).  Recall that a simple random walk $\omega_n$ is typically a distance $O(\sqrt{n})$ from the origin.  There are $O(n^{d/2})$ points near the surface of the $d$-sphere of radius $\sqrt{n}$, and the probability that $\omega_n$ lies at any particular one is $O(n^{-d/2})$:
		$$\max_{z \in \Z^d} \nu(\omega_n = z) = O(n^{-d/2}),$$ 
	where $\nu$ is the simple random walk measure on the space of paths $\Omega$.  This is the $\beta = 0$ case for random polymers; the precise statement for $d \ge 3$ and small $\beta$ is that
		$$\lim_{n\to \oo} \max_{z \in \Z^d} \mu_{n-1}( \omega_n = z ) = 0$$
	for $\PP$-almost every environment $\eta$.  In the strong disorder phase the situation is very different, and polymers are strongly localized.  The same authors proved \cite{carmona2002partition, comets2003directed} that if $d = 1$ or $2$ and $\beta \ne 0$, or if $d \ge 3$ and $\beta$ is large enough, then there exists non-random $c > 0$ such that 
		$$\limsup_{n\to \oo} \max_{z \in \Z^d} \mu_{n-1}( \omega_n = z ) \ge c,$$
	for $\PP$-almost every environment $\eta$.  Giacomin and Toninelli \cite{giacomin2006smoothing} have more recent results on the nature of the phase transition between delocalization and localization.

\subsection{Spin Glasses} \label{spinglasses}

	Spin-glasses are models of interacting particles on a lattice, governed by a Hamiltonian of the form $H(\sigma) = -\sum J_{ij} \sigma_i \sigma_j$.  Unlike a disordered ferromagnet, where the coupling constants $J_{ij}$ are assumed to be random but positive (so that like-spins attract), for a spin-glass model one assumes that the coupling constants $J_{ij}$ can take positive or negative values.  This introduces magnetic frustration (nearby spins need not align), which makes the model difficult to study.
	
	Edwards and Anderson \cite{edwards1975theory} introduced a particularly simple spin-glass model to describe.  Consider a large box $\Lambda \subseteq \Z^d$ of size $|\Lambda| = N^d$, and the space $\Sigma = \{-1,+1\}^\Lambda$ of up-down configurations $\sigma = \{\sigma_i\}$ on $\Lambda$.  Let $J_{ij}$ be an i.i.d. family of random variables on a probability space $(\Omega, \F, \PP)$, and consider the random Hamiltonian $H: \Sigma \to \R$
		$$H_N(\sigma) = - \sum_{|i - j|=1} J_{ij} \sigma_i \sigma_j,$$
	where the sum is over nearest neighbors of $\Lambda$.  For $\beta \ge 0$, consider the random Gibbs measure 
			$$\mu_N(\sigma) = Z^{-1}_N \E^{-\beta H_N(\sigma)}.$$
	with partition function $Z_N = \sum_{\sigma} \E^{-\beta H_N(\sigma)}$.  If the model is ferromagnetic (non-negative coupling constants $J_{ij} \ge 0$), then there are only two ground states:  all sites equal to $+1$ or to $-1$.  In the general spin-glass model, the all-up and all-down states are $\PP$-almost surely no longer ground states.  
	It is an open question if spin-glass models have any non-constant ground states.  As mentioned in Section \ref{gadf}, the two-dimensional disordered ferromagnet ($J_{ij} \ge 0$) is the dual model to Standard FPP, where the interface boundaries of non-trivial ground states \cite{newman1997tds} are the two-sided geodesics on the dual lattice.  Consequently, 	
	
	The spatial structure makes the Edwards-Anderson spin glass extremely difficult to work with.  A drastic simplification is to consider a mean-field model, where the underlying graph is the complete graph $\Lambda$ on $N$ vertices, and every node on the graph interacts with every other one.  Mean-field models are often easier to work with than ones with finite-dimensional interactions.	

	The Sherrington-Kirkpatrick spin-glass \cite{sherrington1975solvable} is one famous example of a mean-field spin-glass model which has exact solutions.  Here the random coupling constants $J_{ij}$ are i.i.d., symmetric random variables with mean zero and variance $J^2 \gg 1$ (when the variance is small, the system is in a weak-disorder phase hence easier to study).  The random Hamiltonian is
		$$H_N(\sigma) = - \frac{1}{\sqrt{N}} \Sigma_{i,j \in \Lambda} J_{ij} \sigma_i \sigma_j,$$
	and the random Gibbs measure 
		$$\mu_N(\sigma) = Z^{-1}_N \E^{-H_N(\sigma)}$$
	with partition function $Z_N = \Sigma_{\sigma} \E^{-H_N(\sigma)}.$  The random free energy is
		$$F_N = - \log Z_N.$$

	Using the non-rigorous technique of replica symmetry breaking, Parisi \cite{parisi1979infinite} calculated an exact form for the free energy in the infinite-volume limit.  Aizenman, Lebowitz and Ruelle \cite{aizenman1987some} rigorously calculated the average value of the free energy per site, as well as the fluctuations.  Talagrand \cite{talagrand1998sherrington} used his powerful concentration-of-measure techniques to rigorously verify Parisi's full ansatz, and Guerra and Toninelli \cite{guerra2002thermodynamic} have pushed these techniques further.

	As discussed in Section \ref{gadf}, it is believed that there are no two-sided minimizing geodesics in two-dimensional Standard FPP.  When interpreted in the context of the disordered ferromagnet, this means that there are no non-trivial ground states, only the unique ground state.  However, this is not believed to be preserved under small perturbations of the metric (a phenomenon called disorder chaos).  It is believed that spin-glass models demonstrate the multiple valleys phenomenon:  there are many very different states which are almost ground states.  Chatterjee \cite{chatterjee2008chaos} has shown that chaos and multiple valleys often occur in tandem in general, along with a phenomenon called ``superconcentration,'' when the variance of the free energy is sublinear.  In \cite{chatterjee2009disorder} Chatterjee proves that the Sherrington-Kirkpatrick model exhibits superconcentration, chaos and multiple valleys.  
	
	Different from the replica method, Mezard, Parisi and Virasoro \cite{mezard1987spin} introduced the ultrametricity assumption to calculate the free energy; see \cite{parisi2000origin} and \cite{aizenman2007mean} for more details.  Derrida \cite{derrida1985generalization} developed the Random Energy Model (REM) to formulate a general proof, which Ruelle \cite{ruelle1987mathematical} improved to the probability cascade technique.  Arguin and Aizenman \cite{arguin2009structure} have recently developed a theory based on multiple valleys to confirm the ultrametricity assumption.
	
	Superconcentration is reminiscent of the sublinear passage-time variation of Benjamini, Kalai and Schramm \cite{benjamini2003first} for Standard FPP discussed in Section \ref{KPZ}, though Chatterjee's demonstration of the phenomenon in the Sherrington-Kirkpatrick model uses very different techniques.  The equivalence of the three phenomena has not yet been shown for Standard FPP.


\section{Included Papers} \label{includedpapers}


	As part of my dissertation work, my advisor Jan Wehr and I wrote the article \cite{lagatta2009shape}, included as Appendix A.  It will be published in the May 2010 issue of the Journal of Mathematical Physics.  This was a collaborative effort between Prof. Wehr and myself.  We wrote this article in order to introduce our continuum model of Riemannian FPP, and demonstrate that we could adapt the basic techniques of Standard FPP for our setting.  The main result is a shape theorem:  large balls in the Riemannian metric grow roughly like Euclidean ones.  As a consequence, we show that the random metric is geodesically complete with probability one.

	I wrote the article \cite{lagatta2010continuous}, included as Appendix B, in order to deal with a conditional probability estimate stemming from the project on geodesics described in Chapter \ref{presentstudy}.  A disintegration (or regular conditional probability) is a way to condition a probability measure on a single point $y$.  In that paper, we introduce continuous disintegrations as those which vary continuously in $y$.  I present a necessary and sufficient condition for continuous disintegrations to exist for Gaussian measures on separable Banach spaces, and analyze how they transform under absolutely-continuous changes of measure.  This project was motivated by the application to Riemannian FPP detailed in Chapter \ref{unifprobproof}; however, the full study of continuous disintegrations was interesting and general enough to warrant submission as a separate publication.
	
	[n.b.: This version of the dissertation does not include the two referenced papers \cite{lagatta2009shape} and \cite{lagatta2010continuous}.]

\chapter{Present Study} \label{presentstudy}

\include{diss_lmgeodesics}

\chapter{Proof of Theorem \ref{frontier}} \label{frontierproof}

\input{diss_proofoffrontierlemma}

\chapter{Sketch of Proof of Claim \ref{unifprob}} \label{unifprobproof}

\input{diss_proofofmarkov}

\bibliographystyle{alpha}
\bibliography{biblio}

\end{document}

%% file: diss_lmgeodesics.tex
	\section{Geometry Background and Notation} \label{geombg}
	
		Before introducing any probabilistic structure, we introduce some geometric notation.  Consider $\R^d$ with $d \ge 2$ and the standard Euclidean coordinates, and fix $\alpha > 0$.  Write
		$$\SPD = \{ \mbox{symmetric, positive-definite $d \times d$ real matrices} \},$$
	and let $g \in C^{2+\alpha}(\R^d,\SPD)$ be a $C^{2+\alpha}$-smooth function on $\R^d$ with values in the space of symmetric, positive-definite matrices.  $g$ defines a Riemannian structure on $\R^d$:  for tangent vectors $v, v' \in T_x \R^d$, we consider the inner product $\langle v, g(x) v' \rangle$.  For a single tangent vector $v$, we denote by $\|v\| = \sqrt{\langle v, g(x) v \rangle}$ and $|v| = \sqrt{\langle v, v \rangle}$ the Riemannian and Euclidean lengths of $v$, respectively.  For a $C^1$-curve $\gamma : [a,b] \to \R^d$, we define the Riemannian and Euclidean lengths of $\gamma$ by
		$$R(\gamma) = \int_a^b \| \dot \gamma(t) \| \sD t \qquad \mathrm{and} \qquad L(\gamma) = \int_a^b | \dot \gamma(t) | \sD t,$$
	respectively.  We say that a curve is finite if it has finite Euclidean length; for our model, Theorem \ref{shapecor} will imply that finite curves have finite Riemannian length.  The Riemannian distance between two points $x$ and $y$ is defined by
		$$d(x,y) = \inf_\gamma R(\gamma),$$
	where the infimum is over all $C^1$-curves $\gamma$ connecting $x$ to $y$.

	A $C^2$-curve $\gamma$ is called a geodesic if it locally minimizes the Riemannian energy functional $R^2(\gamma) = \tfrac{1}{2} \int \|\dot\gamma\|^2$.  Equivalently, geodesics are the solutions to the Euler-Lagrange equation for $R^2$,
		$$\ddot \gamma^k = -\Gamma_{ij}^k \dot \gamma^i \dot \gamma^j,$$
	where we follow the Einstein convention of summing over the repeated indices $i$ and $j$, and where $\Gamma_{ij}^k$ are the Christoffel symbols \cite{lee1997rmi} for the metric $g_{ij}$.  As this is a second-order system of ordinary differential equations, a geodesic is uniquely determined by its starting point and velocity.  Geodesics are locally length-minimizing \cite{lee1997rmi}.  We call a geodesic $\gamma$ minimizing (or globally minimizing) if for all $x, y \in \gamma$, the Riemannian distance $d(x,y)$ is realized as the Riemannian length of $\gamma$ from $x$ to $y$.  Not all geodesics are minimizing; for example, on the sphere, the geodesics are great circles, which are not minimizing past antipodal points.  Geodesics have constant speed \cite{lee1997rmi}; henceforth, we assume $\| \dot \gamma \| = 1$ so that geodesics are parametrized by Riemannian arc length.

	For a Riemannian metric $g$, we define the real, positive functions
		$$\Lambda(x) = \mbox{maximum eigenvalue of $g(x)$} \qquad \mathrm{and} \qquad \lambda(x) = \mbox{minimum eigenvalue of $g(x)$}.$$
	For any $K \subseteq \R^d$, define
		$$\Lambda(K) = \sup_{x\in K} \Lambda(x) \qquad \mathrm{and} \qquad \lambda(K) = \inf_{x\in K} \lambda(x).$$
	By the continuity and positivity of $g$, if $K$ is bounded then
		$$0 < \lambda(K) \le \Lambda(K) < \oo.$$
	For $z \in \Z^d$, let $C_z = [z - 1/2, z + 1/2)^d$ be the unit cube centered at $z$.  Write 
		$$\Lambda_z = \Lambda(C_z) \qquad \mathrm{and} \qquad \lambda_z = \lambda(C_z).$$	


	\section{Riemannian FPP}

		We consider the probability space $\Omega = C^{2+\alpha}(\R^d, \SPD)$ with the $\sigma$-algebra $\F$ generated by cylinder sets.  This space $\Omega$ is a topological subspace of the Fr\'echet space $\hat \Omega = C^{2+\alpha}(\R^d, \Sym)$, where $\Sym$ is the space of symmetric $d\times d$ real matrices with matrix norm.
		
		We call $\Omega$ the space of Riemannian metrics on $\R^d$.  Let $g$ be an $\Omega$-valued random variable with a Radon probability distribution $\PP$ satisfying the following four assumptions:
		\begin{env_ass} \label{assumptions}
		~
		
			\begin{enumerate}[a.]
				\item \label{assumption_stationary} $\PP$ is isotropic, that is, invariant under the isometries of $\R^d$, rotations, translations and reflections.
				\item \label{assumption_finrange} $\PP$ has finite-range dependence.  i.e., there exists $\xi > 0$ such that if $|x-y| \ge \xi$, then $g(x)$ and $g(y)$ are independent.
				\item \label{assumption_finmom} The random variables $\Lambda_0$ and $\Lambda_0/\lambda_0$ have finite moment-generating functions.  That is, 
						$$\EE[\E^{r\Lambda_0}] < \oo \qquad \mathrm{and} \qquad \EE[\E^{r\Lambda_0/\lambda_0}] < \oo $$
				for some $r > 0$.  Since $\lambda_0 \le \Lambda_0$, it follows that $\EE[\E^{r\lambda_0}] < \oo$.
				\item \label{assumption_gaussian} There exists a stationary, mean-zero Gaussian measure $\hat \PP$ on $\hat\Omega$ such that $\PP$ is absolutely continuous with respect to $\hat\PP$ on $\hat \Omega$; and the Radon-Nikodym derivative $\tfrac{\D \PP}{\D \hat \PP}(g)$ is bounded and continuous, and is positive exactly on the open set $\Omega \subseteq \hat\Omega$.
				
			\end{enumerate}
		\end{env_ass}
			
	The first three assumptions imply that $(\Omega, \F, \PP)$ satisfies the hypotheses of \cite{lagatta2009shape}, including a shape theorem with limiting shape equal to a Euclidean ball and almost-sure completeness of the metric $g$.  We summarize these results in this theorem:
		\begin{env_thm} \label{shapecor}
		~
			\begin{enumerate}[a)]
				\item There exists $\mu > 0$ such that $\tfrac{1}{t} d(0,tv) \to \mu$ a.s. and in $L^1$, uniformly in the direction $v \in S^{d-1}$.  Precisely, for all $\epsilon > 0$, with probability one, there exists $T > 0$ such that if $t \ge T$, then $|d(0,tv) - \mu t| \le \epsilon t$ for any $v \in S^{d-1}$.
				\item Let $A = \{ x : |x| \le \mu^{-1} \}$ and $B_t = \{ x : d(0,x) \le t\}$ be the Euclidean and Riemannian balls centered at the origin of radius $\mu^{-1}$ and $t$, respectively.  For all $\epsilon > 0$, with probability one, there exists $T > 0$ such that if $t \ge T$, then
					$$(1-\epsilon) A \subseteq \tfrac{1}{t} B_t \subseteq (1+\epsilon) A.$$
				The Euclidean ball $A$ is called the \emph{limiting shape} of the model.
				\item With probability one, the Riemannian metric $g$ is geodesically complete.  Consequently, with probability one, for all $x$ and $y$ in $\R^d$, there is a finite, minimizing geodesic $\gamma$ connecting $x$ to $y$ such that $d(x,y) = R(\gamma)$.
			\end{enumerate}		
		\end{env_thm}
		\begin{proof}
			The constant $\mu$ is independent of the direction $v$ since the measure $\PP$ is rotationally-invariant.  Part (a) is Proposition 3.3 of \cite{lagatta2009shape}.  Part (b) is Theorem 3.1 of \cite{lagatta2009shape}.  Part (c) is Corollary 3.5 of \cite{lagatta2009shape}.
		\end{proof}


	\section{The evolution of the environment under the geodesic flow} \label{flow}

	For each $x \in \R^d$, the matrix $g(x)$ is positive-definite, hence invertible.  Omitting the $x$, we write $g$ in coordinates as $g_{ij}$, and its inverse $g^{-1}$ as $g^{ij}$.  We define the Christoffel symbols \cite{lee1997rmi}
		$$\Gamma_{ij}^k = \tfrac{1}{2} g^{km} \left( \tfrac{\partial}{\partial x^i} g_{mj} + \tfrac{\partial}{\partial x^j} g_{im} - \tfrac{\partial}{\partial x^m} g_{ij} \right),$$
	where we follow the Einstein notation by summing over repeated indices.  Geodesics are the solutions to the equation 
		$$\ddot \gamma^k = -\Gamma_{ij}^k \dot\gamma^i \dot\gamma^j.$$
	In terms of a vector field $U : \R^{2d} \to \R^{2d}$, geodesics are the flow lines for
		$$U(x,v) = (v, -\Gamma_{ij}^k(x) v^i v^j \E_k),$$
	where $\E_k$ is the $k^{\mathrm{th}}$ standard basis vector in $\R^d$.  Let $F_t : \R \times \R^{2d} \to \R^{2d}$ be the geodesic flow, so that $\tfrac{\D}{\D t} F_t = U(F_t)$ and
		$$F_t(x,v) = (\gamma_{x,v}(t), \dot\gamma_{x,v}(t)),$$
	where $\gamma_{x,v}$ is the unique geodesic starting at $x$ in the direction $v$.  By assumption, the metric $g(x)$ is $C^{2+\alpha}$-smooth, so the Christoffel symbols $\Gamma_{ij}^k$ and the field $U$ are $C^{1+\alpha}$-smooth.  Consequently, the flow $F_t$ is $C^{1+\alpha}$-smooth \cite{arnold1988geometrical} hence locally Lipschitz.

	Fix $v \in S^{d-1}$.  Rather than fixing the environment $g \in \Omega$ and considering the flow $F_t(0,v)$ along the geodesic $\gamma_v := \gamma_{0,v}$ (the Eulerian perspective), we instead consider a reference frame centered along a particle traveling along $\gamma_{v}$ (the Lagrangian perspective).  Define the random flow $\sigma_{t} : \Omega \to \Omega$ on the space of Riemannian metrics by
		$$(\sigma_{t} g)(u) = g(u + \gamma_{v}(t)).$$
	The variable $u \in \R^d$ represents the displacement from $\gamma_v(t)$, so that $(\sigma_t g)(0) = g(\gamma_v(t))$ always represents the metric at $\gamma_v(t)$.  The flow $\sigma_t$ induces the random measure $\PP \circ \sigma_t^{-1}$ on $\Omega$.

	\begin{env_cla} \label{abscondt}
		Fix $v \in S^{d-1}$.  With probability one, for all $t \in \R$ the random measure $\PP \circ \sigma_{t}^{-1}$ on $\Omega$ is absolutely continuous with respect to $\PP$.  That is, there exists a family of measurable functions $\rho_t : \Omega \to \R$ so that for all measurable $f : \Omega \to \R$,
			$$\int_\Omega f(\sigma_t g) \sD \PP(g) = \int_\Omega f(g) \rho_t(g) \sD \PP(g).$$
	\end{env_cla}
	
	This should follow from the main theorem of Geman and Horowitz \cite{geman1975random}; see Zirbel \cite{zirbel2001lagrangian} for a more recent presentation.  They call a vector field homogeneous if its law is translation-invariant.  By the isotropy of $\PP$, the vector field $U(x,v)$ is homogeneous in the first coordinate.  If $U$ were homogeneous in both coordinates, then the claim would immediately follow by Proposition 8.2 of \cite{zirbel2001lagrangian}.  However, this is not the case, and their work must be modified for this situation.




	\section{Rarity of Minimizing Geodesics}
	
	Consider the set of all minimizing, unit-speed geodesics between the origin and the boundary sphere of the Euclidean ball $\Beuc_n := \Beuc(0,n)$ of radius $n$.  Let $\V_n \subseteq S^{d-1}$ be the set of initial velocities of these geodesics.  Note that these geodesics may exit $\Beuc_n$.  Clearly, $\V_n$ is monotonically decreasing in the sense that $\V_{n+1} \subseteq \V_n$.  Let $\V = \bigcap \V_n$.

	\begin{env_pro} \label{Vnonempty}
		The set $\V$ is non-empty and closed.
	\end{env_pro}
	\begin{proof}
		Since the metric is complete with probability one, each distance $d(0, n\E_1)$ is realized by a finite minimizing geodesic $\gamma_n$ connecting $0$ to $n\E_1$.  Let $v_n \in S^{d-1}$ be the initial velocities of these geodesics.  Since the sphere is compact, a subsequence $v_{n_k}$ converges to some $v \in S^{d-1}$.  Let $\gamma$ be the unique geodesic with $\gamma(0) = 0$ and $\dot\gamma(0) = v$, parametrized by Riemannian length.  We claim that $\gamma$ is minimizing.
		
		Let $x = \gamma(t)$ and $x' = \gamma(t')$ be two points along the curve $\gamma$.  As mentioned in the previous section, the geodesic flow is smooth in the initial conditions, so
	 		$$x = \lim_{k\to\oo} \gamma_{n_k}(t) \qquad \mathrm{and} \qquad x' = \lim_{k\to\oo} \gamma_{n_k}(t').$$
	 	Since the distance function $d$ is continuous and the finite geodesics $\gamma_{n_k}$ are minimizing,
	 		$$d(x,x') = \lim_{k\to\oo} d(\gamma_{n_k}(t), \gamma_{n_k}(t')) = |t-t'|.$$
	 	This proves that $\gamma$ globally minimizes length.
	 	
	 	In fact, the same argument shows that $\V$ is closed.  Let $v_n \in \V$, and suppose that $v_n \to v$ in $S^{d-1}$.  Let $\gamma_n$ and $\gamma$ be the geodesics starting at the origin in directions $v_n$ and $v$, respectively.  The above argument shows that $\gamma$ is minimizing, so $v \in \V$.
	\end{proof}
	
	We call $w$ an asymptotic direction of $\gamma$ if the limit of $\gamma(t)/|\gamma(t)|$ exists and equals $w$ as $t \to \oo$.  
	Howard and Newman \cite{howard1997euclidean} have shown that for their rotationally-invariant model of Euclidean first-passage percolation, every one-sided minimizing geodesic has an asymptotic direction.  The key to their proof is that the limiting shape of Euclidean FPP is a Euclidean ball.  Since the global structure of both their model and our own are similar, we conjecture that the same is true in this setting:


	\begin{env_con} \label{asympconj}
		For every $v \in \V$, there exists $w \in S^{d-1}$ such that
			$$\lim_{t\to\oo} \frac{\gamma_v(t)}{|\gamma_v(t)|} = w.$$

%
	\end{env_con}

	

	If this conjecture holds, then we should be able to improve non-emptyness of $\V$ to uncountability:
	\begin{env_cor} \label{uncountableV}
		If Conjecture \ref{asympconj} holds, then the set $\V$ is uncountable.
	\end{env_cor}
	
	We claim that for each $w \in S^{d-1}$, there is some $v \in \V$ so that $\gamma_v$ has asymptotic direction $w$.  The proof mimics that of Proposition \ref{Vnonempty}:  we begin with the minimizing geodesics $\gamma_n$ from $0$ to $nw$, and take the limit of a subsequence of their starting directions $\dot \gamma_{n_k}(0) \to v$.  By Conjecture \ref{asympconj}, $\gamma_v$ has an asymptotic direction $w'$.  An additional argument is needed to show that $w = w'$; for example a result that the transversal fluctuation exponent $\xi$ is less than 1, as for lattice FPP models \cite{howard2004mfp}.\newline

%
%
%


	The main goal of this project is to show that for a deterministic $v \in S^{d-1}$, the geodesic $\gamma_v$ starting in direction $v$ is length minimizing with probability zero.  
	
	\begin{env_cla} \label{mainresult}
		For each $v \in S^{d-1}$, the event 
			\begin{equation} \label{vmin}
				\{ v \in \V\} = \{\mbox{$\gamma_v$ is minimizing}\} \end{equation}
		has probability zero.
	\end{env_cla}
	
	This is a new result for this model for which there is no analogue in lattice FPP models.  At the time of submitting this dissertation, we do not yet have a full proof of this claim.  However, we have formulated the general argument of the proof, as well as many technical lemmas.  We plan to work through all the technical details and submit this soon for publication.

	In Section \ref{sect_frontiertimes}, we prove that minimizing geodesics are transient, i.e. leave every compact set.  We use this along with some results on dependent lattice FPP which we developed in \cite{lagatta2009shape} to prove a global existence statement:  there are a number of ``frontier times'' along a minimizing geodesic $\gamma_v$ at which things are ``well-behaved.''  In Section \ref{sect_unifprob}, using the continuous disintegrations which we developed for stationary Gaussian measures in \cite{lagatta2010continuous}, we show that at each of these times $t_k$, there is a uniform probability estimate for a destabilizing phenomenon to occur in front of $\gamma_v(t_k)$.  In Section \ref{sect_bump}, we argue that this phenomenon should be a ``bump surface,'' to exploit positive curvature so that the curve $\gamma_v$ develops conjugate points.
	
	By a standard application of Tonelli's theorem \cite{folland1999real}, Claim \ref{mainresult} can be improved to show that, with probability one, $\V$ has measure zero on the sphere $S^{d-1}$.  This is not just a technicality:  in light of Corollary \ref{uncountableV}, we believe that $\V$ is uncountable.  We further believe that $\V$ has no isolated points and is nowhere dense, so that the set of minimizing directions is a random topological Cantor set.

	\begin{env_thm} \label{Vmeas0}
	 	If Claim \ref{mainresult} holds, then with probability one, the set $\V$ has measure zero on the sphere $S^{d-1}$.  Precisely, if $\nu$ is the uniform measure on $S^{d-1}$, then
	 		$$\PP\big( \nu(\V) = 0 \big) = 1.$$
	\end{env_thm}
	
	\begin{proof}
		For $v \in S^{d-1}$, let $E_v = \{v \notin \V\}$ be the event that the geodesic $\gamma_v$ is not minimizing.  Claim \ref{mainresult} implies that $\PP(E_v) = 1$.  Write $\V^c = \{ v \in S^{d-1} : E_v ~\mbox{occurs} \}$ for the directions which do not give minimizing geodesics, and let $\nu$ be the uniform measure on $S^{d-1}$.  Tonelli's theorem \cite{folland1999real} implies that
			\begin{eqnarray*}
				\int_\Omega \nu(\V^c) \sD \PP(\omega)
				&=& \int_\Omega \nu(v : E_v ~\mbox{occurs} ) \sD \PP(\omega) = \int_\Omega \int_{S^{d-1}} 1_{E_v}(\omega) \sD \nu(v) \sD \PP(\omega) \\
				&=& \int_{S^{d-1}} \int_\Omega 1_{E_v}(\omega) \sD \PP(\omega) \sD \nu(v) = \int_{S^{d-1}} \PP(E_v) \sD \nu(v) = \int_{S^{d-1}} 1 \sD \nu(v) = 1,
			\end{eqnarray*}
		since $\PP(E_v) = 1$. Thus $\nu(\V^c) = 1$ with probability one, so $\nu(\V) = 0$.
	\end{proof}

	

	\section{Transience of Geodesics and Existence of Frontier Times} \label{sect_frontiertimes}
	

	As part of their definition in Section \ref{geombg}, geodesics are parametrized by Riemannian arc length, so $\|\dot \gamma(t)\| = 1$ for all $t$.  This is the natural parametrization from the point of view of differential geometry, as it depends only on the intrinsic geometry.  In our probabilistic model, the initial Euclidean coordinate system is also natural.  Since geodesics are curves in $\R^d$, we will also consider them parametrized by Euclidean arc length $l$, so that $|\dot \gamma(l)| = 1$ for all $l$.

	The following theorem demonstrates that minimizing geodesics are transient, whether parametrized by Riemannian or Euclidean length.  We show that for any (possibly random) compact set $K$, there exists a uniform time after which all minimizing geodesics never return to $K$.  We use the notation $\gamma_v$ to mean the unique geodesic starting at $0$ in direction $v \in S^{d-1}$.
	
	\begin{env_thm} \label{transientgeodesics}
	~
	\begin{enumerate}[a)]
		\item Suppose that geodesics are parametrized by Riemannian arc length $t$.  With probability one, if $K$ is a (possibly random) compact set in $\R^d$, then there exists a time $T$ such that for all $v \in \V$ and $t > T$, $\gamma_v(t) \notin K$.  
		
		\item Suppose that geodesics are parametrized by Euclidean  arc length $l$.  With probability one, if $K$ is a (possibly random) compact set in $\R^d$, then there exists a time $L$ such that for all $v \in \V$ and $l > L$, $\gamma_v(l) \notin K$.  
	\end{enumerate}
	\end{env_thm}

	We require almost-sure completeness of the metric in our proof of part (b), where we assume that a Riemannian ball of finite radius must be compact in $\R^d$.
	


	
	\begin{proof}[Proof of a)]
		Let $\hat K = \Beuc(0,r)$ be the smallest Euclidean ball centered at the origin which contains $K$.  The metric $g$ is continuous hence bounded on the ball $\hat K$, so the maximum eigenvalue $\Lambda(\hat K)$ is finite.  Let $T = r \sqrt{\Lambda(\hat K)}$.
		
		Let $v \in \V$ and suppose that $\gamma_v$ is the unique geodesic starting at the origin in direction $v$.  If $\gamma_v(t) \in K$ for some time $t$, then since $\gamma_v$ is minimizing,
			$$t = d(0, \gamma_v(t)) \le r\sqrt{\Lambda(\hat K)} = T,$$
		where we estimate the distance by the Riemannian length of the straight-line path between $0$ and $\gamma_v(t)$.  Thus, if $t > T$, then $\gamma_v(t) \notin K$.  
	\end{proof}
		
	\begin{proof}[Proof of b)]
		For a (possibly random) compact set $K$, let $T$ be as in part (a).  Consider $B = \Briem(0,T)$, the closed Riemannian ball centered at the origin of radius $T$.  By the almost-sure completeness of the metric, $B$ is compact.  The metric is positive-definite and continuous, so $\lambda(B) > 0$.  Let $L = T / \sqrt{\lambda(B)}$.
		
		Let $v \in \V$ and suppose that $\gamma_v(l) \in K$; we will show that $l \le L$.  Let $t(l) = \int_0^l \sqrt{\langle \dot \gamma_v, g \dot \gamma_v \rangle}$ be the Riemannian arc length of $\gamma_v$ from $0$ to $l$.  Since $\gamma_v(l) \in K$, the above argument shows that $t(l) \le T$.  Furthermore, since $t$ is an increasing function of $l$, for all $l' \le l$ the Riemannian times $t(l')$ are bounded above by $T$, hence $\gamma_v(l') \in B$.  Thus
			$$T \ge t(l)  = \int_0^l \sqrt{\langle \dot \gamma_v, g \dot \gamma_v \rangle} \ge l \sqrt{\lambda(B)},$$
		since $\gamma_v$ is parametrized by Euclidean arc length so $\langle \dot\gamma_v, \dot\gamma_v \rangle = 1$. Therefore, $l \le T / \sqrt{\lambda(B)} = L$.
	\end{proof}
	
	The next theorem is an improvement on the previous one.  Not only are minimizing geodesics transient, but for each $v \in \V$, there exists a sequence of ``frontier times'' $t_k(v) \uparrow \oo$ such that things are ``nice'' at $\gamma_v(t_k)$.  First, the geodesic satisfies a cone condition at these times:  there is a uniform $\theta < \tfrac{\pi}{2}$ such that the angle between $\gamma_v(t_k)$ and $\dot \gamma_v(t_k)$ is less than $\theta$.  In particular, this means that at $t_k$, the geodesic is not tangent to the Euclidean ball of radius $|\gamma_v(t_k)|$ centered at the origin.  Second, there is a uniform upper bound on the $C^{2+\alpha}$-norm of the metric $g$ in a uniform neighborhood $B_k$ of $\gamma_v(t_k)$, as well as a lower bound on $\lambda$, the minimum eigenvalue of $g$.
	
	Let $\xi$ be the finite-dependence length of the metric.  i.e., if $|u - v| \ge \xi$, then $g(u)$ and $g(v)$ are independent.

		\begin{env_thm} \label{frontier}
			There exist non-random $\beta \in (0, 1)$ and $h > 0$ such that, with probability one, for all $v \in \V$, there exists a sequence of ``frontier times'' $t_k(v) \uparrow \oo$ such that
				\begin{itemize}
					\item The angle between $\gamma_v(t_k)$ and $\dot \gamma_v(t_k)$ is at most $\theta := \cos^{-1} \beta$, uniformly in $k$.  i.e., 
						$$\langle \gamma_v(t_k), \dot \gamma_v(t_k) \rangle \ge |\gamma_v(t_k)| \, |\dot\gamma_v(t_k)| \, \cos \theta = |\gamma_v(t_k)| \, |\dot \gamma_v(t_k)| \, \beta.$$
					\item Let $\rho = 2\xi / \beta$, for $\xi$ as above.  Write
						\begin{equation} \label{Bk}
							B_k = \Beuc(\gamma(t_k), \rho) \end{equation}
					for the Euclidean ball of radius $\rho$ centered at $\gamma(t_k)$.  Then
						$$\|g\|_{C^{2+\alpha}(B_k)} + \tfrac{1}{\lambda(B_k)} \le h.$$
				\end{itemize}
		\end{env_thm}

		The event in this theorem holds simultaneously for all directions in the set $\V$ with probability one, though the particular sequence of times $t_k$ depends on the direction $v$.  In fact, we will prove this theorem for arbitrary $\rho$ in \eqref{Bk}, though in that case the non-random constant $h$ will depend on $\rho$.  The proof is technical, and uses some lemmas from \cite{lagatta2009shape}.  It can be found in Appendix \ref{frontierproof}.

	\section{Uniform Probability Estimates at Frontiers} \label{sect_unifprob}
	

	For this section, we fix $v \in S^{d-1}$, and consider the unique geodesic $\gamma_v$ starting from the origin in direction $v$.  If $\gamma_v$ is to be minimizing, a necessary condition will be that there is a sequence of ``frontier times'' along the geodesic.  We argue in Claim \ref{unifprob} that at each of these times, there is a uniform probability $p$ with which a certain event occurs. 
	
	To this more precise, we consider a filtration, ordered by space rather than time.  Since minimizing geodesics are transient, a natural filtration to consider is
		$$\F_r := \sigma \{ g(x) : |x| \le r \},$$
	the $\sigma$-algebra generated by the metric in the closed Euclidean ball $\Beuc_r = \Beuc(0,r)$.

	Define the random function $\tau_v : [0,\oo) \to [0,\oo]$ as the first time that $\gamma_v$ leaves the ball of radius $r$.  That is,
		$$\tau_v(r) = \inf \left\{ t : |\gamma_v(t)| = r ~\mathrm{and~is~increasing} \right\},$$
	where $\tau_v(r) = \oo$ if $\gamma_v$ is trapped in the ball $\Beuc_r$ for all time (i.e. $|\gamma_v(t)| \le r$ for all $t$).  Where they are finite, the random functions $\tau_v(r)$ are all strictly increasing and right-continuous with left limits.  The exit times $\tau_v(r)$ depend only on the metric in the Euclidean ball of radius $r$, hence are adapted to the filtration $\F_r$.  For transient geodesics $\gamma_v$ (including minimizing geodesics by Proposition \ref{transientgeodesics}), the exit time $\tau_v(r)$ is finite for all $r$.
	
	Let $\beta$, $h$ and $\rho = 2\xi/\beta$ be as in Theorem \ref{frontier}, and fix $v \in S^{d-1}$.  We will call $R$ a frontier of $\gamma_v$ if the exit time $t := \tau_v(R)$ is finite and satisfies the conclusions of Theorem \ref{frontier}, where $B_k$ is replaced by 
		$$B = \Beuc(\gamma(t), \rho) \cap \Beuc_R,$$
	the part of the neighborhood around $\gamma(t)$ which is contained in the large ball $\Beuc_R = \Beuc(0,R)$. 
	
	\begin{env_def}
		We define $R\ge 0$ to be a frontier of $\gamma_v$ if the exit time $t := \tau_v(R)$ is finite, the angle between $\gamma_v(t_k)$ and $\dot \gamma_v(t_k)$ is at most $\theta := \cos^{-1} \beta$, uniformly in $k$, and
			\begin{equation} \label{hestimate}
				\|g\|_{C^{2+\alpha}(B)} + \tfrac{1}{\lambda(B)} \le h. \end{equation}
	\end{env_def}

	Frontiers are ``stopping times'' (in the probabilistic sense) with respect to the filtration $\F_r$, since the event $$\{\mbox{$R$ is a frontier of $\gamma_v$}\} \cap \{R \le r\}$$ depends only on the metric in the ball $\Beuc_r$ (i.e. the event is $\F_r$-measurable).  Theorem \ref{frontier} implies that there is a sequence of frontiers along minimizing geodesics:
	
	\begin{env_cor} \label{Vfrontiers}
		With probability one, if $v \in \V$, then there is a sequence of frontiers $R_k \uparrow \oo$ along $\gamma_v$.
	\end{env_cor}
	
	Let $v \in S^{d-1}$.  We will use frontiers to test if $v \in \V$.  If we can not find a sequence of frontiers $R_k$ along $\gamma_v$, then Corollary \ref{Vfrontiers} implies that $v \not \in V$.  If there does exist such a sequence $R_k$, then our Claim \ref{unifprob} will imply that there is a uniform probability $p$ so that at each frontier time $R_k$, the geodesic $\gamma_v$ encounters a phenomenon which destabilizes the minimization property.\newline

	Let $\mathcal O_x : \R^d \to \R^d$ be a family of affine transformations on $\R^d$ which map $0 \mapsto x$ and $-|x|\E_1 \mapsto 0$.\footnote{For example, let $\mathcal O_x^{\mathrm{trans}}$ be the translation which sends $0$ to $|x|\E_1$, and let $\mathcal O_x^{\mathrm{rot}}$ be the identity transformation if $x$ is parallel to $\E_1$; otherwise, let $\mathcal O_x^{\mathrm{rot}}$ be the rotation which fixes the $(d-2)$-dimensional space $\span\{\E_1, x\}^\perp$, and rotates the vector $\E_1$ in the plane $\span\{\E_1,x\}$ to be parallel to $x$.  Define $\mathcal O_x = \mathcal O_x^{\mathrm{rot}} \mathcal O_x^{\mathrm{trans}}$.}  
	Fix $v \in S^{d-1}$, and define the $\F_r$-measurable affine transformation $\mathcal O_r := \mathcal O_{\gamma_v(\tau_v(r))}$ on the event $\{\tau_v(r) < \oo\}$.  The map $\mathcal O_r$ rotates and translates $\R^d$ so that at the frontier time $\tau_v(r)$, the transformed geodesic is sitting at the origin with the former ball $\Beuc(0,r)$ contained entirely in the left half-space.  We define the random transformation $\mathcal O_r$ on the space $\Omega = C^{2+\alpha}(\R^d,\SPD)$ by 
		$$(\mathcal O_r g)(u) := g (\mathcal O_r u ), \qquad u \in \R^d.$$
	If we consider a particle traveling along the geodesic $\gamma_v$, then by adopting the point of view of the particle, $\Theta_r g$ is the environment the particle sees at time $\tau_v(r)$.  The left half-space represents the ``past'' of the particle's trajectory, and the right half-space the ``future.''  The transformation $\O_r$ is a random shift, followed by a random rotation.  Consequently, the random measure $\PP \circ \O_r^{-1}$ on $\Omega$ is absolutely continuous with respect to $\PP$, as in Section \ref{flow}.
	
	\begin{env_cla} \label{abscondr}
		Fix $v \in S^{d-1}$.  With probability one, for all $r \ge 0$ the random measure $\PP \circ \O_r^{-1}$ on $\Omega$ is absolutely continuous with respect to $\PP$.  That is, there exists a family of measurable functions $\rho_r : \Omega \to \R$ so that for all measurable $f : \Omega \to \R$,
			\begin{equation} \label{chmeasshift}
				\int_\Omega f(\O_r g) \sD \PP(g) = \int_\Omega f(g) \rho_r(g) \sD \PP(g). \end{equation}
	\end{env_cla}
	
	This should follow from Claim \ref{abscondt}, where we must account for the stopping time $\tau_v(r)$, as well as the random rotation.\newline

	For $\rho = 2\xi/\beta$ as above, let $B_\oo = \Beuc(0,\rho) \cap \{x : x^1 \le 0 \}$ be the closed left half-ball of radius $\rho$, and let 
		\begin{equation} \label{Bdef}
			B_r = \Beuc(0,\rho) \cap \Beuc(-r\E_1, r) \end{equation}
	be the part cut out of $B_\oo$ by the large ball $\Beuc(-r\E_1,r)$.
	Fix $\eta > 0$, and define the cone $C$ in the right half-space by
		\begin{equation} \label{Cdef}
			C = \left\{ x \in \R^d : 0 \le x^1 \le \eta \mathrm{~and~} \sqrt{(x^2)^2 + \dots + (x^d)^2} \le \rho x^1 \right\}. \end{equation}
	If $\phi$ denotes the angle of of the cone $C$ from the horizontal axis, then $\cos \phi = \beta/2$.  Thus $\phi$ is strictly greater than $\theta = \cos^{-1} \beta$, since cosine is decreasing.
	
	Write $W = B_\oo \cup C$.  Note that the only points in the left half-space which are Euclidean distance less than $\xi$ away from $C$ are those in $B_\oo$.  Conditioned on the left half-space, the metric $g|_C$ in the cone depends only on the metric $g|_{B_\oo}$ in the half-ball.  This is an important point which we exploit in the proof of Claim \ref{markov} to show that there is a Markov Property of the metric at frontier times.
	
	Define $\Theta_r : \Omega \to C^{2+\alpha}(W, \SPD)$ by
		\begin{equation} \label{Thetar}
			(\Theta_r g)(u) = g\left( \mathcal O_r u \right), \qquad u \in W. \end{equation}
	Thus $\Theta_r g$ is the metric in the neighborhood of $\gamma_v(\tau_v(r))$, rotated and translated to lie at the origin.  Let
		$$\eta_r : C^{2+\alpha}(W, \SPD) \to C^{2}(B_r, \SPD)$$		
	be the restriction-and-inclusion map, defined by $(\eta_r x)(u) = x(u)$ for $u \in B_r$.  The map $\eta_r \Theta_r : \Omega \to U_r$ is $\F_r$-measurable.

	\begin{env_cla} \label{unifprob}

		Let $v \in S^{d-1}$, and let $\{\Theta_r\}$ be the family of $\F_r$-adapted random maps as defined in \eqref{Thetar}.  If $U \subseteq C^{2+\alpha}(W, \SPD)$ is open, then there exist non-random $p > 0$ and $r_0 > 0$ such that if $R \ge r_0$ is a frontier of $\gamma_v$ and $\eta_R \Theta_R g \in \eta_R U$, then
			$$\PP \left(\Theta_{R}^{-1} U \ |\  \F_{R} \right) > p.$$
	\end{env_cla}
	
	
	
	

	The event $\{\eta_R \Theta_R g \in \eta_R U\}$ is simply that the part of the metric $g$ contained in $\Beuc_R$ is compatible with the event $\Theta^{-1}_R U$.  In the sequel, this event will be implied by the estimate \eqref{hestimate}.  
	
	We sketch the proof of this claim, which involves some tools coming from probability in Banach spaces and developed in \cite{lagatta2010continuous}.  Assumption \ref{assumptions}.\ref{assumption_gaussian} of this model was that $\PP$ is absolutely continuous with respect to a Gaussian measure, which implies that the disintegration (i.e. regular conditional probability) satisfies certain continuity properties \cite[Theorem 11]{lagatta2010continuous}.  The Arzel\`a-Ascoli theorem \cite{folland1999real} implies that the set of metrics for which \eqref{hestimate} holds is compact in the $C^2$-norm.  This gives us a positive lower bound for the event to occur.  The proof is technical and can be found in Appendix \ref{unifprobproof}.
	

\section{Construction of a Bump Surface at Frontier Times and Proof of Main Result} \label{sect_bump}

	Consider the cone $C$ as defined in \eqref{Cdef} as a manifold with boundary.  Let $Z = C^2(C, \SPD)$ be the space of $C^2$-Riemannian metrics on $C$.  Let $\phi$ be the angle of the cone at $0$, so that $\tan \phi = \rho/\xi = 2/\beta$, and $\phi$ is strictly greater than $\theta = \cos^{-1} \beta$.  Consequently, if a geodesic $\gamma$ starts at the origin with initial rightward direction within angle $\theta$ of $\E_1$, that is,
		$$\dot \gamma^1(0) \ge |\dot \gamma(0)| \cos \theta = |\dot \gamma(0)| \beta,$$
	then $\gamma(t)$ is in the interior of the cone for small, positive time $t$.
	
	Let $Y = C^2(B_\oo, \SPD)$ be the space of $C^2$-Riemannian metrics on the half-ball $B_\oo$, defined in \eqref{Bdef} (this is the space $Y_\oo$ as defined in Appendix \ref{unifprobproof}).  The set $\Gamma = \{g \in Y : \|g\|_{C^{2+\alpha}} + 1/\lambda \le h \}$ is compact in $Y$ by the Arzel\`a-Ascoli theorem \cite{folland1999real}.

	\begin{env_cla} \label{bumpclaim}
		There exists a continuous map $b : \Gamma \to X$ and $\epsilon > 0$ such that if $\| g - b(g|_{B_\oo})\|_Z < \epsilon$, then for all geodesics $\gamma$ starting at 0 with initial directions within an angle $\theta$ of $\E_1$, there exists a point $x$ in the interior of $C$ such that $0$ and $x$ are conjugate points along $\gamma$.
	\end{env_cla}
	
	For each $g \in \Gamma$, the function $b(g) : C \to \SPD$ is a Riemannian metric on the cone $C$, which we call a ``bump metric.''  All the geodesics which pass over the bump develop conjugate points \cite{lee1997rmi} and lose the minimization property.  While we will see this exact Riemannian manifold with probability zero, the loss of minimization persists under small perturbations of the metric.

	This construction has two elements:  first that we can construct a Riemannian metric $\tilde g := b(g)$ such that the geodesics remain in the cone $C$ and develop conjugate points, and that this is stable under a uniform perturbation $\epsilon$ of the metric.  We have not yet completed the construction with all the technical details, but we include the sketch of our argument here.  The cone $C$ meets the half-ball $B_\oo$ at the origin, so the Riemannian metric $\tilde g$ must agree at $0$ with $g$ up to second derivatives.  Since $\Gamma$ is compact, these derivatives are all bounded.  Other than this condition, we have absolute freedom to choose a Riemannian metric which does whatever we want in $C$.
	
	Let $\tilde \Gamma_{ij}^k$ be the Christoffel symbols \cite{lee1997rmi} for the Riemannian metric $\tilde g$, so that the geodesic equation is
		$$\ddot \gamma^k = -\tilde \Gamma_{ij}^k \dot \gamma^i \dot \gamma^j,$$
	where we follow the Einstein summation convention and sum over the repeated indices $i$ and $j$.  In particular, for the first coordinate
		$$\ddot \gamma^1 = -\tilde \Gamma_{ij}^1 \dot \gamma^i \dot \gamma^j.$$
	As a geodesic approaches the boundary of the cone, we want it to be accelerated rightward, so we want the Christoffel symbols $\tilde \Gamma_{ij}^1$ to be negative and very large near the boundary.  Once we guarantee that the geodesics are moving roughly parallel and to the right, we smooth the metric out into a spherical metric.  This is the origin of the name ``bump'':  the attached Riemannian manifold begins with arbitrary (but bounded) positive, zero or negative curvature at the origin, then as geodesics follow the manifold the curvature becomes constant and positive.  It is the presence of positive curvature which forces geodesics to develop conjugate points, after which they are not minimizing \cite{lee1997rmi}.
	
	Conjugate points occur when the solution to the Jacobi equation \cite{lee1997rmi} along a geodesic vanishes twice.  The Jacobi equation is a differential equation with coefficients comprised of the second derivatives of the metric $\tilde g$.  Consequently, zeros to solutions are stable under small $C^2$ perturbations of the metric.  For each $y \in \Gamma$, let $\epsilon(y) > 0$ be the maximum such perturbation such that the consequence of Claim \ref{bumpclaim} holds.  This should be a continuous function of $y$ in the compact set $\Gamma$, hence the minimum $\epsilon = \inf_{y \in \Gamma} \epsilon(y)0$ is non-zero.

	\begin{env_lem}
		The set $U \subseteq X$ defined by
			$$U = \{ g \in X : \| g|_{C} - b(\eta_\oo g) \|_Z < \epsilon \}$$
		is open in $X$.
	\end{env_lem}
	\begin{proof}
		The function $f : X \to \R$ defined by
			$$f(g) = \| g|_{C} - b(\eta_\oo g) \|_{Z}$$
		is continuous, and $U = f^{-1}( (-\oo, \epsilon) )$.
	\end{proof}
	
	Finally, we can prove the main result of the paper, and show that $v \in \V$ with probability zero.
	
	\begin{proof}[Proof of Claim \ref{mainresult}]
	
		Let $v \in S^{d-1}$.  If there is no sequence of frontiers $R_k \uparrow \oo$ along $\gamma_v$, then $v \notin \V$ by Corollary \ref{Vfrontiers}.  Suppose that the event does hold, and let $R_k$ be the sequence of frontiers.  Let $U$ be as in the preceding lemma, so if any of the events $\Theta_{R_k}^{-1} U$ occur then the geodesic $\gamma_v$ is not minimizing.  
	
	For all $k$,
		$$\PP \left( \bigcap_{k'=1}^k \left(\Theta_{R_{k'}}^{-1} U \right)^c ~\Big|~ \F_{R_k} \right) \le (1-p)^k.$$
	
	Thus with probability one, the event $\Theta_{R_k}^{-1} U$ occurs for some $k$.

	\end{proof}

%% file: diss_proofoffrontierlemma.tex

	In this Appendix we prove Theorem \ref{frontier}.  Corollary \ref{shapepaperlemmas} is a summary of some results from \cite{lagatta2009shape}.  We apply those results in the proof of Lemma \ref{euclengthboundlemma}, which controls the Euclidean arc length of a minimizing geodesic.  The key assumption is that $\Lambda/\lambda$---the ratio of the largest eigenvalue of the Riemannian metric $g$ in a unit cube to the smallest eigenvalue---is a random variable with strong tail decay properties.  This means that for most cubes it passes through, a minimizing geodesic will not wiggle too much.

	We recall some notation from \cite{lagatta2009shape}.  	For $z \in \Z^d$, we write $z = (z^1, \dots, z^d)$.  We say that $z, z' \in \Z^d$ are $*$-adjacent if $\max_{1\le i\le d} (z - z')^i \le 1$.  The $*$-lattice is the graph with vertex set $\Z^d$, and edge set given by $*$-adjacency; that is, the usual lattice $\Z^d$ along with all the diagonal edges.
		
	We say that a set $\Gamma \subseteq \Z^d$ is $*$-connected if for all $z, z' \in \Gamma$, there is a path from $z$ to $z'$ along the $*$-lattice which remains in the set $\Gamma$.  Technically, that there is a finite sequence of $*$-adjacent points beginning with $z$ and ending with $z'$, all contained in $\Gamma$.  
	
	Let $X_z$ be a stationary, non-negative random field on the $*$-lattice with finite-range dependence, and with a finite moment-generating function
			\begin{equation} \label{finmom2}
				M(r) = \EE[\E^{r X}] < \oo \qquad \mbox{for all $r \in \R$.} \end{equation}
	The finite-range dependence means that there exists $\xi > 0$ such that if $|z - z'| \ge \xi$, then $X_z$ and $X_{z'}$ are independent.  We write
		$$X(\Gamma) = \sum_{z\in\Gamma} X_z.$$

	Assumption \ref{assumption_finmom} implies that $\Lambda$ and $\Lambda/\lambda$ have finite moment-generating functions and satisfy \eqref{finmom2}.  Since $\lambda < \Lambda$, $\lambda$ also satisfies \eqref{finmom2}.

	\begin{env_cor} \label{shapepaperlemmas}
	~
	\begin{enumerate}[a)]
		\item For $\mu$ as in Corollary \ref{shapecor}, with probability one, there exists $M_1 > 0$ such that if $|x| \ge M_1$, then $d(0,x) \le 2\mu|x|$.
	
		\item Suppose that $X_z$ is stationary and positive, and satisfies finite-range dependence and \eqref{finmom2}.  For any $A > 0$ there is a non-random $B > 0$ such that, with probability one, there exists $N > 0$ such that for all $n \ge N$, if $\Gamma$ is a $*$-connected set containing the origin and $X(\Gamma) \le An$, then $|\Gamma| \le Bn.$
	
		\item Suppose that $X_z$ is stationary and non-negative, and satisfies finite-range dependence and \eqref{finmom2}.  For any $B > 0$ there is a non-random $C > 0$ such that, with probability one, there exists $N > 0$ such that for all $n \ge N$, if $\Gamma$ is a $*$-connected set containing the origin and $|\Gamma| \le Bn$, then $X(\Gamma) \le Cn.$
	\end{enumerate}
	\end{env_cor}
	\begin{proof}
		Part (a) is implied by Theorem \ref{shapecor}.a.  Parts (b) and (c) are Lemmas 2.2 and 2.3 of \cite{lagatta2009shape}, respectively, applied to the constant sequence $a_n \equiv 0$.
	\end{proof}

	\begin{env_lem} \label{euclengthboundlemma}
		There exists a non-random $D \ge 1$ such that, with probability one, there exists $M > 0$ such that if $|x| \ge M$ and $\gamma$ is a length-minimizing geodesic connecting $0$ to $x$, then
			\begin{equation} \label{euclengthbound}
				|x| \le L(\gamma) \le D|x|, \end{equation}
		where $L(\gamma)$ denotes the Euclidean length of $\gamma$ between $0$ and $x$.
	\end{env_lem}
	
	\begin{proof}
		The lower estimate $|x| \le L(\gamma)$ is trivial, since $\gamma$ has Euclidean length at least that of the straight line path from $0$ to $x$.

		By Corollary \ref{shapepaperlemmas}.a, with probability one, there exists $M_1 > 0$ such that if $|x| > M_1$, then
			$$d(0,x) \le 2\mu |x|.$$
		
		Apply Corollary \ref{shapepaperlemmas}.b to $A = 8\mu$ and $X_z = \lambda_z$.  Thus there exists a non-random $B > 0$ such that, with probability one, there exists $N_1 > 0$ such that for all $n \ge N_1$, if $\Gamma$ is a finite $*$-connected set which contains the origin and $\lambda(\Gamma) \le 8\mu n$, then $|\Gamma| \le Bn$.
		
		By Assumption \ref{assumption_finmom}, $\Lambda_z/\lambda_z$ has a finite moment-generating function.  Apply Corollary \ref{shapepaperlemmas}.c to the above $B$ and $X_z = \Lambda_z/\lambda_z$.  Thus there exists a non-random $C > 0$ such that, with probability one, there exists $N_2 > 0$ such that for all $n \ge N_2$, if $\Gamma$ is a finite $*$-connected set which contains the origin and $|\Gamma| \le Bn$, then $(\Lambda/\lambda)(\Gamma) \le Cn$.
		
		Set $D = \tfrac{3^d B}{2} + 2C\sqrt{d}$, and let $|x| \ge \max\{M_1, N_1, N_2,1\}$.  Let $n$ be the smallest integer greater than $|x|$; we will later use the trivial estimate $n \le 2|x|$.  Let $\gamma$ be a length-minimizing geodesic between $0$ and $x$.  Since $\gamma$ connects the origin to a point Euclidean distance $|x|$ away, $L(\gamma) \ge |x|$.  Define the discrete set
			\begin{equation} \label{Gammadef}
				\Gamma = \{ z \in \Z^d : L(\gamma \cap C_z) \ge 1/4 \}. \end{equation}
		That is, $z \in \Gamma$ if $\gamma$ spends at least Euclidean length $1/4$ in the unit cube $C_z$.  The set $\Gamma$ is $*$-connected; see the discussion following (2.8) of \cite{lagatta2009shape}.  Clearly, $0 \in \Gamma$.
		
		Since $\gamma$ is length-minimizing, 
			$$R(\gamma) = d(0,x) \le 2\mu|x| \le 2\mu n.$$
		Furthermore, by summing $\lambda_z$ over the points of $\Gamma$, we get an upper bound using $R(\gamma)$:
			$$\tfrac{1}{4} \lambda(\Gamma) \le \sum_{z \in \Gamma} L(\gamma \cap C_z) \lambda_z \le \sum_{z \in \Gamma} R(\gamma \cap C_z) \le R(\gamma) \le 2\mu n.$$
		Thus, $\lambda(\Gamma) \le 8\mu n$, hence 
			\begin{equation} \label{GammaB}
				|\Gamma| \le Bn, \end{equation}
		and $(\Lambda/\lambda)(\Gamma) \le Cn$.
		
		In each cube $C_z$, we can estimate the Euclidean length of $\gamma$ using $\Lambda_z/\lambda_z$:
			$$L(\gamma \cap C_z) \lambda_z \le R(\gamma \cap C_z) \le \Lambda_z \sqrt{d},$$
		so
			$$\qquad L(\gamma \cap C_z) \le \frac{\Lambda_z}{\lambda_z} \sqrt{d}.$$

		Define the set $\hat \Gamma$ consisting of $\Gamma$ and all neighboring points on the $*$-lattice:
			$$\hat \Gamma = \{ z \in \Z^d : \exists ~ z' \in \Gamma \mbox{~s.t. $z$ and $z'$ are $*$-adjacent} \} \supset \Gamma.$$
		The geodesic $\gamma$ is completely contained in union of the cubes with centers $\hat \Gamma$.  The geodesic can get contributions to Euclidean length from the cubes with centers $z \in \hat \Gamma \backslash \Gamma$, but only up to $1/4$ and there are fewer than $3^d |\Gamma| \le 3^d Bn$ of such cubes.  Thus
			$$L(\gamma) \le \sum_{z \in \hat\Gamma \backslash \Gamma} L(\gamma \cap C_z) + \sum_{z \in \Gamma} L(\gamma \cap C_z) \le \frac{3^d B}{4} n + \sqrt{d} \sum_{z \in \Gamma} \frac{\Lambda_z}{\lambda_z} \le \frac{3^d B}{4} n + \sqrt{d} Cn = \tfrac{1}{2} Dn,$$
		since $D = \tfrac{3^d B}{2} + 2C\sqrt{d}$.  Since $n \le 2|x|$, the proof is complete.
	\end{proof}
	
	Let $\beta = 1/2D < 1$, and let $\theta \in (0, \pi/2)$ be the angle such that $\cos \theta = \beta$.  For $v \in \V$, consider the length-minimizing geodesic $\gamma_v$, and suppose that it is parametrized by Euclidean arc length $l$.  Write $r_v(l) = |\gamma_v(l)|$.  Define the set of Euclidean frontier times of $\gamma_v$ to be
		$$F_v = \left\{l : \dot r_v(l) > \beta \mathrm{~and~} r_v(l) = \sup_{l' \le l} r_v(l') \right\}.$$
	
	In the next lemma, we show that the set of Euclidean frontier times takes up a non-zero fraction of the Euclidean length of $\gamma_v$.

	\begin{env_lem} \label{frontierdensity}
		With probability one, for all $v \in \V$, the set of Euclidean frontier times $F_v \subseteq [0,\oo)$ comprises right-open intervals and is unbounded.  Furthermore, there exists non-random $\delta > 0$ such that, with probability one, there exists $L > 0$ such that if $l \ge L$, then
			$$\Leb(F_v \cap [0,l]) \ge \delta l$$
		for all $v \in \V$.
	\end{env_lem}
	\begin{proof}

		We first argue that $F_v$ is right-open.  Suppose $l \in F_v$.  Since $\dot r$ is continuous, there exists $\epsilon > 0$ such that if $h \in [0,\epsilon)$, then $\dot r_v(l + h) > \beta$.  Since $r_v(l) = \sup_{l' \le l} r_v(l')$ and $r$ is strictly increasing on $[l, l+\epsilon)$, $r_v(l+h)$ is the new supremum.  Thus $[l,l+\epsilon) \subseteq F_v$.
	
		Let $D$ and $M$ be as in Lemma \ref{euclengthboundlemma}.  Let $K = \Beuc(0,M)$ be the Euclidean ball of (random) radius $M$.  By Theorem \ref{transientgeodesics}, with probability one, all minimizing geodesics escape $K$ in uniform time:  there exists $L$ such that if $l\ge L$ and $v \in \V$, then $r_v(l) = |\gamma_v(l)| \ge M$, hence
			$$r_v(l) \le l \le D r_v(l).$$
		
		Let $\delta = 1/(2D - 1)$.  Write $S = \left\{l : r_v(l) = \sup_{l' \le l} r_v(l') \right\}$ for the times $l$ at which $r_v(l)$ attains the supremum, so that we can decompose the non-frontier times $F_v^c$ by
			$$F_v^c = \left( \{0 \le \dot r \le \beta \} \cap S \right) \cup S^c.$$
		If $l \ge L$, then the fundamental theorem of calculus implies that
			\begin{equation} \label{3ints}
			D^{-1} l \le r_v(l) = \int_0^l \dot r = \int_{F_v \cap [0,l]} \dot r + \int_{\{0 \le \dot r \le \beta\} \cap S \cap [0,l] } \dot r + \int_{S^c \cap [0,l]} \dot r. \end{equation}

		
		Since $f_v(l) := \sup_{l' \le l} r_v(l') - r_v(l)$ is continuous, $S^c = f_v^{-1}((0,\oo))$ is open, hence a union of open intervals.  Let $I$ be a maximal subinterval of $S^c$.  The curve $\gamma_v$ is transient by Theorem \ref{transientgeodesics} so $f(l) = 0$ for arbitrarily large $l$; this implies that $I$ is bounded.  At both endpoints of $I$, the function $r$ equals $\sup r$, so the third integral of \eqref{3ints} vanishes.  
		
		Write $b(l) = \Leb(F_v \cap [0,l])$; we must show $b(l) \ge \delta l$.  Since the geodesic is parametrized by Euclidean length, $\dot r \le 1$.  We use this to estimate the first integral of \eqref{3ints}; for the second integral, we use $\dot r \le \beta$.  Thus
			$$D^{-1} l \le 1 \cdot b(l) + \beta \cdot (l - b(l)) + 0.$$
		Since $\beta = 1/2D$ and $\delta = 1/(2D-1)$, by rewriting this expression, we have $b(l) \ge \delta l$ as desired.
	\end{proof}

%
	

	\begin{proof}[Proof of Theorem \ref{frontier}]
		Suppose that geodesics are parametrized by Euclidean length.  For $v \in S^{d-1}$, let $t_v(l)$ be the change in parametrization to Riemannian arc length along $\gamma_v$, and let $r_v(l) = |\gamma_v(l)|$.  It suffices to prove that for all $v \in \V$, there exists a sequence $l_k \uparrow \oo$ such that the conclusions of Theorem \ref{frontier} hold for the sequence $t_k := t_v(l_k)$.  The metric is complete with probability one, so $l_k \uparrow \oo$ implies that $t_k \uparrow \oo$.

		Let $v \in \V$, and fix $l \in F_v$.  We first prove that the angle between $\gamma_v(l)$ and $\dot \gamma_v(l)$ is less than $\theta := \cos^{-1} \beta$.  This follows quickly from the definition of frontier times and elementary trigonometry.  Since $\gamma_v$ is parametrized by Euclidean length, $|\dot \gamma_v(l)| = 1$.  Since $l$ is a frontier time, $\dot r_v(l) \ge \beta$: the projection of $\dot \gamma(l)$ onto the direction $\gamma(l)$ is at least $\beta$.  Consequently, the angle between $\dot\gamma(l)$ and $\gamma(l)$ is at most $\theta$, where $\cos \theta = \beta$.

		Fix $\rho > 0$, and let $m$ be the minimum number of cubes $C_z$ which can cover any Euclidean ball of radius $\rho$.  Let $B$ be as in Lemma \ref{euclengthboundlemma}.  Write $\tilde\rho = \rho + \sqrt{d}$.

		Let $v \in \V$.  Define an increasing sequence of frontier times $l_j \in F_v$ and balls $B_j \subseteq \R^d$ as follows.  Let $l_0 = 0$ and
			$$l_j = \inf \left\{ l \in F_v : l > l_{j-1} \mathrm{~and~} |\gamma(l) - \gamma(l_{j'})| \ge 2\tilde\rho \mathrm{~for~} j' < j \right\}.$$
		Define the ball $B_j = \Beuc(\gamma(l_j), \rho)$ of radius $\rho$ centered at $\gamma(l_j)$, and let
			$$\tilde B_j = \{ z \in \Z^d : B_j \cap C_z \ne \emptyset\}$$
		be the centers of the cubes $C_z$ which form a discrete cover of $B_j$, so $|\tilde B_j| \le m$.  The discrete sets $\tilde B_j$ are disjoint, since two distinct $\rho$-balls $B_j$ are separated by distance at least $\sqrt{d}$.


		\begin{env_lem}
			Let $\delta$ and $L$ be as in Lemma \ref{frontierdensity}, and let $A = \tfrac{8\tilde\rho}{\beta \delta}$.  If $l_j \ge L$ then
				\begin{equation} \label{ljj}
					l_j \le Aj. \end{equation}
		\end{env_lem}
		\begin{proof}
			Clearly, the balls $B(\gamma(l_j), 4\tilde \rho)$ of larger radius $4\tilde\rho$ cover the image under $\gamma$ of all frontier times $F_v$:
				$$F_v \subseteq \bigcup_{j'=1}^\oo \left\{ l \in F_v : |\gamma(l) - \gamma(l_{j'})| \le 4\tilde\rho \right\} =: \bigcup_{j'=1}^\oo I_{j'},$$
			hence
				$$\Leb(F_v \cap [0,l_j]) \le \sum_{j'=1}^j \Leb(I_{j'}).$$
			On $I_j$, the maximum distance to the origin $\sup_{l' \le l} r_v(l')$ can grow by at most $8\tilde\rho$, the diameter of the ball $\Beuc(\gamma(l_j), 4\tilde\rho)$.  Thus by the fundamental theorem of calculus,
				$$8 \tilde\rho \ge \int_{I_j} \dot r \ge \beta \, \Leb(I_j).$$
			If $l_j \ge L$, then Lemma \ref{frontierdensity} implies that $\delta l_j \le \Leb(F_v \cap [0,l_j])$.  Thus
				$$l_j \le \tfrac{1}{\delta} \sum_{j'=1}^j \Leb(I_{j'}) \le \tfrac{8\tilde\rho}{\delta \beta} j = Aj.$$
		\end{proof}
		
		Let 
			$$W_j = \left\{y \in \R^d : |y - \gamma(l)| \le \rho \mathrm{~for~some~} l \in [0,l_j] \right\}$$
		be the $\rho$-neighborhood of $\gamma|_{[0,l_j]}$.  Let 
			$$\tilde\Gamma_j = \{ z \in \Z^d : C_z \cap W_j \ne \emptyset\}$$
		be the centers of the cubes $C_z$ which cover $W_j$.  Note that $\tilde B_{j'} \subseteq \tilde\Gamma_j$ for all $j' \le j$.

		\begin{env_lem} \label{tildeGammaj}
			There exists non-random $B' > 0$ and there exists $J_1 > 0$ such that if $j \ge J_1$, then
				$$|\tilde\Gamma_j| \le B' j.$$ 					
		\end{env_lem}
		\begin{proof}
			Let $L$ be as in Lemma \ref{frontierdensity}.  Let			
				$$\Gamma_j = \{ z \in \Z^d : \gamma|_{[0,l_j]} \cap C_z \ne \emptyset \}$$
			be the centers of the cubes $C_z$ which the curve $\gamma|_{[0,l_j]}$ meets.  As in Lemma \ref{euclengthboundlemma}, there exists non-random $B > 0$ and there exists $L_1 > 0$ such that if $l_j \ge L_1$, then $|\Gamma_j| \le B l_j$.  Let $J_1$ be the minimum $j$ such that $l_j \ge \max\{L,L_1\}$, and suppose $j \ge J_1$.  By \eqref{ljj}, $l_j \le Aj$.  Let $B' = mBA$, so
				$$|\tilde\Gamma_j| \le m |\Gamma_j| \le mB l_j \le mBA j = B'j.$$

		\end{proof}

		Now let $h \in (0,\oo)$, and let $A_z^h$ be the event that
			\begin{equation} \label{badcube}
				\|g\|_{C^{2+\alpha}(C_z)} + \tfrac{1}{\lambda(C_z)} > h; \end{equation}
		Let $X_z^h = 1(A_z^h)$ be the indicator function of the event $A_z^h$.  Since the family $X_z^h$ only takes the values $0$ and $1$, it is bounded hence has a finite moment-generating function.  
		
		Apply Corollary \ref{shapepaperlemmas}.c to the $B'$ from Lemma \ref{tildeGammaj} and the family $X_z^h$.  Thus there exists a non-random $C(h) > 0$ (depending on $h$) such that, with probability one, there exists $J_2 > 0$ such that for all $j \ge J_2$, if $\Gamma$ is a finite $*$-connected set which contains the origin and $|\Gamma| \le B'j$, then $X^h(\Gamma) \le C(h) j$.  
		
		With probability one, the metric $g$ is $C^{2+\alpha}$ and positive everywhere.  Thus for every $z \in \Z^d$, 
			$$\lim_{h \to \oo} \PP(A_z^h) = 0.$$
		Consequently, $C(h) \to 0$ as $h\to\oo$.  Choose a value of $h$ large enough so that
			$$C(h) < \frac 1 2.$$
		
		Let $j \ge \max\{J_1, J_2\}$.  By the above lemma, $|\tilde\Gamma_j| \le B' j$ so
			$$X^h(\tilde\Gamma_j) \le C(h) j < \frac{j}{2}.$$
		That is, the number of points $z \in \tilde\Gamma_j$ for which that the event $A_z^h$ occurs is fewer than $j/2$.  There are $j$ disjoint sets $\{\tilde B_{j'}\}$ contained in $\tilde\Gamma_j$; consequently, there are at least $j/2$ balls $B_{j_k}$ such that
			$$\| g \|_{C^{2+\alpha}(B_{j_k})} + \tfrac{1}{\lambda(B_{j_k})} \le h.$$
		As $j \to \oo$, we may choose infinitely many $j_k \to \oo$.  This proves Theorem \ref{frontier}. 
	\end{proof}

%% file: diss_proofofmarkov.tex

	Let $B_r$ be the half-ball and $C$ be the cone defined as in \eqref{Bdef} and \eqref{Cdef}, respectively, and let $W = B_\oo \cup C$.  Recall that $\Omega = C^{2+\alpha}(\R^d, \SPD)$.  Write $\Sym$ for the space of symmetric $d \times d$ real matrices, and consider the Banach spaces
		$$X = C^{2+\alpha}(W, \Sym) \qquad \mathrm{and} \qquad Y_r = C^{2}(B_r, \Sym), \quad r \le \oo$$
	equipped with the $C^{2+\alpha}$ and $C^2$ norms, respectively.  The set inclusions $B_r \subseteq B_\oo \subseteq W \subseteq \R^d$ induce restriction-and-inclusion maps
		$$\xymatrix{ \Omega \ar[r]^\chi & X \ar[r]^{\eta_\oo} \ar[rd]_{\eta_r} & Y_\oo \ar[d]^{\varphi_r} \\ && Y_r }$$	
	We need to account for the parameter $r$ in our maps, since the region we will be conditioning on later on will be cut out from the large ball $\Beuc(-r\E_1, r)$.

	Let $\PP_X = \PP \circ \chi^{-1}$ be the push-forward of the probability measure $\PP$ on $X$.  Similarly, let $\PP_{Y_r} = \PP_X \circ \eta_r^{-1}$ be the push-forward probability measures on $Y_r$ for $r \le \oo$.  When there is no ambiguity we will write $\PP$ for $\PP_X$.  These measures satisfy the change-of-variable equations
		\begin{equation} \label{chvarboth}
			\int_\Omega f(\chi g) \sD \PP(g) = \int_X f(x) \sD \PP_X(x) \qquad \mathrm{and} \qquad \int_X g(\eta_r x) \sD \PP_X(x) = \int_{Y_r} g(y) \sD \PP_{Y_r}(y) \end{equation}
	for any measurable functions $f : X \to \R$ and $g : Y_r \to \R$.  
		
%

	Assumption \ref{assumptions}.\ref{assumption_gaussian} implies that there exists a stationary, mean-zero Gaussian measure $\hat \PP$ on $X$ such that $\PP$ is absolutely continuous with respect to $\hat\PP$ on $X$; and the Radon-Nikodym derivative 
		$$\rho(x) := \tfrac{\D \PP}{\D \hat \PP}(x)$$
	is bounded and continuous, and is positive exactly on the open subset 
		$$X^0 := C^{2+\alpha}(W,\SPD) \subseteq X.$$
	This implies that $X^0$ has full $\PP$-measure.  Let $ \PP_{Y_r} = \PP \circ \eta_r^{-1}$ denote the push-forward measures of $\PP$ on the spaces $Y_r$.  Consequently, the sets 
		$$Y_r^0 := C^{2+\alpha}(B_r, \SPD) = \eta_r(X^0) \subseteq Y_r$$
	have full $\PP_{Y_r}$-measure, though not open since the spaces $Y_r$ are equipped with the $C^2$-norm instead of the $C^{2+\alpha}$-norm.

	\begin{env_pro} \label{contdis}
		There exist regular conditional probabilities $\nu_r : Y_r^0 \times \B(X) \to [0,1]$ such that:
		\begin{enumerate}[a)]
			\item If $\Gamma \subseteq Y_r^0$ is compact in $Y_r$, and if $y_n \in \Gamma$ and $y_n \to y$, then the measures $\nu_r(y_n, \cdot)$ converge weakly to $\nu_r(y, \cdot)$ on $X$.
			
			\item If $B \subseteq X^0$ is open and $y \in Y_r^0 \cap \eta_r(B)$, then $\nu_r(y,B) > 0$.
			
			\item \textbf{Claim:}  If $\Gamma \subseteq Y_\oo^0$ is compact in $Y_\oo$ and $B \subseteq X^0$ is open, then for all $\epsilon > 0$, there exists $R > 0$ such that if $r \ge R$, then for all $y \in \Gamma$,
				\begin{equation} \label{approxmeas}
					\nu_r( \varphi_r y, B) \ge \nu_\oo( y, B) - \epsilon. \end{equation}
		\end{enumerate}
	\end{env_pro}
	
	\begin{proof}[Proof of a)]
		Let $c : W \times W \to \Sym$ be the matrix-valued covariance function of the Gaussian measure $\hat \PP$.  i.e., if $\omega \in X$ is a realization of $\hat \PP$, then
			$$c(u,v) = \hat\EE( \omega(u) \omega(v) ),$$
		where the product is matrix multiplication.  For $u \in W$, write $c_u(\cdot) = c(u,\cdot)$, so $c_u \in X$.  If $u \in B_r$, then $c_u \in Y_r$.

		Let $K : X^* \to X$ be the covariance operator of the Gaussian measure $\hat \PP$ on $X$, defined by $Kf(u) = f(c_u)$.  The spaces 
			$$\hat Y_r = \overline{\eta_r K \eta_r^* Y_r^*} \subseteq Y_r$$  		
		have full $\hat \PP_{Y_r}$-measure.  Furthermore, as subspaces of $X$,
			$$\overline{\bigcup K\eta_r^* Y_r^*} = \overline{ \span_{u \in \bigcup B_r} \{c_u\} } = \overline{ \span_{u \in B_\oo} \{c_u\} }= \overline{ K \eta_\oo^* Y_\oo^* },$$
		since the family $\{c_u\}$ is equicontinuous for $u \in \bigcup B_r$.

		Let $u_0 \in \bigcap B_r$.  We assume that $\hat\PP$ is non-degenerate, so $\|c(u_0,u_0)\|_{\Sym} > 0$.  Let $M = \sup_{u_\in B_\oo} \|c_u\|_X / \|c(u_0,u_0)\|_{\Sym} < \oo$.  Since $c$ is stationary,
			$$M_r := \sup_{e \in Y_r^*} \frac{\|K\eta_r^* e\|}{\|\eta_r K \eta_r^* e\|} = \sup_{u \in B_r} \frac{\|c_u\|_X}{\|c_u\|_Y} \le \frac{\sup_{u_\in B_\oo} \|c_u\|_X}{\|c(u_0,u_0)\|_{\Sym}} = M < \oo$$
		uniformly for all $r \le \oo$.  On the dense subspace $\eta_r K \eta_r^* Y_r^*$ of $\hat Y_r$, define $m_r : \eta_r K \eta_r^* Y_r^* \to X$ by $y \mapsto \eta_r^{-1}(y)$.  This linear map has operator norm $M_r \le M < \oo$, hence we may extend $m_r$ continuously to all of $\hat Y_r$.  
				
		For all $r \le \oo$ and $y \in \hat Y_r$, let $\PP_{Y_r}^y$ be the Gaussian measure on $X$ with mean $m_r(y)$ and covariance operator $\hat K_r = K - K\eta_r^* m_r^*$.  By Theorem 6 of \cite{lagatta2010continuous}, each $\hat \PP^y_{Y_r}$ is a continuous disintegration on $\hat Y_r$ with respect to the map $\eta_r$.  That is, $\hat \PP^y_{Y_r}$ is a regular conditional probability with respect to $\eta_r$, and if $y_n \in \hat Y_r$ and $y_n \to y$, then $\hat \PP^{y_n}_{Y_r} \to \hat\PP^y_{Y_r}$ weakly.
			
		
		In the context of $\PP \ll \hat\PP$, Theorem 11 of \cite{lagatta2010continuous} implies that $\PP_{Y_r} \ll \hat\PP_{Y_r}$ with bounded, continuous Radon-Nikodym derivative 
			$$\rho_{Y_r}(y) := \tfrac{\D \PP_{Y_r}}{\D \hat\PP_{Y_r}}(y) = \int_{\eta_r^{-1}(y)} \rho(x) \sD \hat\PP^y_{r} (x).$$
		The function $\rho_{Y_r}$ is positive exactly on the set $Y_r^0$.  Theorem 11 of \cite{lagatta2010continuous} also implies that the measure $\nu_r : Y_r^0 \times \B(X) \to [0,1]$ defined by
			$$\nu_r(y,B) = \int_{B} \frac{\rho(x)}{\rho_{Y_r}(y)} \sD \hat\PP^y_r(x).$$
		is a regular conditional probability for $\PP$ and that property (a) holds.
	\end{proof}
		
	\begin{proof}[Proof of b)]
		Let $B \subseteq X^0$ be open, and let $y \in Y_r^0 \cap \eta_r(B)$.  Choose $x_0 \in B$ such that $\eta_r(x_0) = y$.  The function $\rho$ is positive and continuous at $x_0$ and the set $B$ is open so there exist $a > 0$ and $\delta > 0$ such that the open ball $B(x_0, \delta)$ is contained in $B$, and $\rho(x) > a$ on $B(x_0, \delta)$.  Gaussian measures assign positive measure to open sets, so
			$$\nu_r(y,B) = \int_{B} \frac{\rho(x)}{\rho_{Y_r}(y)} \sD \hat\PP^y_r(x) \ge \frac{a}{\rho_{Y_r}(y)} \hat\PP^y_r (B(x_0, \delta) ) > 0.$$
	\end{proof}
	
	\begin{proof}[Sketch of proof of c)]
		Define $Y_{r,\oo} = \eta_\oo K \eta_r^* Y_r^* \subseteq \hat Y_\oo$.  The union of these spaces is dense in $\hat Y_\oo$, since
			$$\overline{\bigcup Y_{r,\oo} } = \overline{\span_{u\in \bigcup B_r} \{c_u\} } = \overline{\span_{u\in B_\oo} \{c_u\} }= \hat Y_\oo.$$
		We show now that the operators $m_r \varphi_r$ converge uniformly to $m_\oo$ on $\hat Y_\oo$.  On $Y_{r,\oo}$, the maps $m_r \varphi_r$ and $m_\oo$ are equal, since if $g \in Y_r^*$,
			$$(m_r \varphi_r - m_\oo)(\eta_\oo K \eta_r^* g) = m_r \eta_r K \eta_r^* g - m_\oo \eta_\oo K \eta_\oo^* (\varphi_r^* g) = 0.$$		
		
		Let $\epsilon > 0$, and choose $r_0$ such that for all $y \in \hat Y_\oo$, there is some $y' \in Y_{r_0,\oo}$ such that $\|y - y'\|_{Y_\oo} \le \epsilon / 2M$.  
		Then
			$$\| (m_r \varphi_r - m_\oo)(y) \|_X \le \| (m_r \varphi_r - m_\oo)(y - y') \|_X + 0 \le ( \|m_r\| \|\varphi_r\| + \|m_\oo\| ) \|y-y'\| \le \epsilon,$$
		proving that $m_r \varphi_r$ converges uniformly to $m_\oo$.
		
		The Gaussian measure $\hat \PP^y_{Y_\oo}$ has mean $m_\oo(y)$ and covariance operator $\hat K_\oo = K - K\eta_\oo^* m_\oo^*$, and the Gaussian measure $\hat \PP^{\varphi_r y}_{Y_r}$ has mean $m_r( \varphi_r y)$ and covariance operator
			$$\hat K_r = K - K\eta_r^* m_r^* = K - K (\varphi_r \eta_\oo)^* m_r^* = \hat K_\oo - K \eta_\oo^* (m_r \varphi_r - m_\oo)^*.$$
		It should follow from standard theory on Gaussian measures \cite{bogachev1998gm, ibragimov1978grp} that an approximation statement like \eqref{approxmeas} holds for the Gaussian measures $\hat \PP^y_{Y_\oo}$ and $\hat \PP^{\varphi_r y}_{Y_r}$.  Once we have proved that, estimate \eqref{approxmeas} should follow easily from the explicit construction of $\nu$ from the Gaussian disintegration.
	\end{proof}

	Let $\chi : \Omega \to X$ be the restriction-and-inclusion map from $\Omega$ to $X$ and let $\Theta_r = \chi \mathcal O_r : \Omega \to X$ be the family of $\PP$-random maps as defined in \eqref{Thetar}, depending on fixed $v \in S^{d-1}$.


	\begin{env_pro}[Markov Property] \label{markov}
		If $v \in S^{d-1}$ and $f : X \to \R$ is measurable, then for all $r < \oo$,
			\begin{equation} \label{markoveqn}
				\EE\left(f \circ \Theta_r | \F_r \right) = \int_X f(x) \, \nu_r(\eta_r \Theta_r g, \D x ) \end{equation}
		on the event $\{ \tau_v(r) < \oo \}$ for $\PP$-almost every $g \in \Omega$.
	\end{env_pro}
	\begin{proof}
		For this proof, we suppose that $r$ is fixed, and consequently drop it from our notation when it is clear.  Recall that $\Theta_r = \chi \circ \O_r$.  By Claim \ref{abscondr}, the measure $\PP \circ \O_r^{-1}$ is absolutely continuous to $\PP$ on the event $\{\tau_v(r) < \oo\}$, so we will first prove a statement analogue to \eqref{markoveqn} without the random transformation $\O_r$.  After that, we will transform the measure and prove \eqref{markoveqn}.
		
		Consider the $\sigma$-algebras
		$$\F_{\Beuc} = \sigma\{ g(x) : x \in \Beuc(-r\E_1, r) \}, \quad \F_B = \sigma\{ g(x) : x \in B_r \}, \quad \mathrm{and} \quad \F_W = \sigma\{ g(x) : x \in B_r \cup C \}.$$
		By the construction of $B_r$ and $C$, the sets $\Beuc(-r\E_1, r) \backslash B_r$ and $C$ are separated by Euclidean distance at least $\xi$.  Thus as Hilbert subspaces of $L^2(\Omega, \F)$, this implies that
			$$L^2(\Omega, \F_{\Beuc}) \cap L^2(\Omega, \F_W) = L^2(\Omega, \F_B).$$
		The random variable $f \chi : \Omega \to \R$ is $\F_W$-measurable, so conditioning it on the $\sigma$-algebra $\F_{\Beuc}$ reduces to conditioning on $\F_B$:
			\begin{equation} \label{fchi}
				\EE(f\chi | \F_{\Beuc}) = \EE(f\chi | \F_B). \end{equation}
		
		Now, we claim that 
			\begin{equation} \label{markovconnection}
				\EE(f\chi | \F_B ) = \int_X f(x) \nu(\eta_r \chi g, \D x) \end{equation}
		for $\PP$-almost every $g$.  Suppose $A$ is a $\F_B$-measurable event.  The map $(\eta \chi)^{-1} (\eta \chi) : \F \to \F$ projects an event onto the coordinates generated by points in $U$; consequently, $(\eta \chi)^{-1} (\eta \chi)A = A$.  Thus by applying both change-of-variable formulas \eqref{chvarboth} and the disintegration equation,
			\begin{eqnarray*}
				\int_A \EE(f\chi | \F_B) \sD \PP(g) &=& \int_A f(\chi g) \sD \PP(g) \\
				&=& \int_{\chi A} f(x) \sD \PP_X(x) \\
				&=& \int_{\eta \chi A} \int_X f(x) \, \nu(y, \D x) \sD \PP_Y(y) \\
				&=& \int_A \int_X f(x) \, \nu(\eta \chi g, \D x) \sD \PP(g).
			\end{eqnarray*}
		This proves \eqref{markovconnection}.\newline
		
		We return to the random-transformation case to prove \eqref{markoveqn}.  Let $A \in \F_r$, and write $A' = A \cap \{\tau_v(r) < \oo\}$.  We claim that 
			\begin{equation} \label{markoveqnintegrated}
				\int_{A'} \EE\left(f \circ \Theta_r | \F_r \right) \sD \PP(g) = \int_{A'} \int_X f(x) \, \nu_r(\eta_r \Theta_r g, \D x ) \sD \PP(g). \end{equation}
		The left side is equal to
			\begin{equation} \label{frankie}
				\int_{A'} f(\chi \O_r g) \sD \PP(g) = \int_\Omega f(\chi g) 1_{A'}(\O_r^{-1} g) \rho_r(g) \sD \PP(g) \end{equation}
		by the change of measure \eqref{chmeasshift}.  The random transformation $\O_r^{-1}$ on $\Omega$ is $\F_{\Beuc}$-measurable, as is the function $\rho_r$.  Consequently, the right-hand side of \eqref{frankie} is equal to
			$$\int_\Omega \EE( f\chi \cdot 1_{\O_r A'} \cdot \rho_r | \F_{\Beuc} ) \sD \PP = \int_\Omega \EE( f\chi | \F_{\Beuc}) 1_{\O_r A'} \cdot \rho_r \sD \PP = \int_\Omega \EE( f\chi | \F_B)  1_{\O_r A'} \cdot \rho_r \sD \PP$$
		since $\EE( f\chi | \F_{\Beuc}) = \EE( f\chi | \F_B)$ by \eqref{fchi}.  Substituting \eqref{markovconnection}, this is equal to
			$$\int_\Omega \left(\int_X f(x) \, \nu(\eta_r \chi g, \D x) \right) 1_{\O_r A'}(g) \rho_r(g) \sD \PP(g) = \int_{A'} \int_X f(x) \, \nu(\eta_r \chi \O_r g, \D x) \sD \PP(g),$$
		where we transform the measure back to $\PP$ via \eqref{chmeasshift}.

	\end{proof}

	\begin{env_cla}[Strong Markov Property] \label{strongmarkov}
		Supposing that Claim \ref{markov} holds, if $v \in S^{d-1}$ and $f : X \to \R$ is measurable, then
			$$\EE\left(f \circ \Theta_R | \F_R \right) = \int_X f(x) \, \nu_R(\eta_R \Theta_R g, \D x )$$	
		on the event $\{ \mbox{$R$ is a frontier of $\gamma_v$} \}$ for $\PP$-almost every $g \in \Omega$.
	\end{env_cla}

	This proof follows the classic proof of the Strong Markov Property \cite{durrett1996probability}, where we approximate the random frontiers by deterministic radii.  We have not yet worked through the argument in full detail, but there should be no technical complications. With the Strong Markov Property in hand, we are ready to prove Claim \ref{unifprob}:
	
	\begin{proof}[Proof of Claim \ref{unifprob}]
		Let $v \in S^{d-1}$ and let $B \subseteq C^{2+\alpha}(W, \SPD)$ be open.  Since we are considering frontiers, define
			$$\Gamma = \{ y \in Y : \|y\|_{C^{2+\alpha}(B_\oo)} + \tfrac{1}{\lambda(B_\oo)} \le h \} \subseteq Y_\oo$$
		for the value of $h$ as in Theorem \ref{frontier}.  Because of the H\"older condition $\alpha$ on the second derivatives, the Arzel\`a-Ascoli Theorem \cite{folland1999real} implies that $\Gamma$ is compact in $Y_\oo$.  
		
		
		Let 
			$$p = \tfrac{1}{2} \inf_{y \in \Gamma} \nu_\oo(y, B).$$
		Since $B$ is open, Proposition \ref{contdis}.b implies that the function $\nu_\oo(\cdot, B)$ is lower semi-continuous.  Hence on the compact set $\Gamma$ it attains its minimum $2p$.  By Proposition \ref{contdis}.a, this is positive so $p > 0$.

		By the Strong Markov Property,
			$$\PP ( \Theta_R^{-1} B | \F_r ) = \nu_R (\eta_R \Theta_R g, B )$$
		on the event $\{ \mbox{$R$ is a frontier of $\gamma_v$} \}$.  This event further implies that $\eta_R \Theta_R g \in \varphi_R \Gamma$.

		Following the discussion on the definition of $m_r$ in the proof of Proposition \ref{contdis}.a, define the continuous map $\alpha_r = \varphi_r^{-1} : Y_r \to Y_\oo$ on the dense subspace $\varphi_r \eta_\oo K \eta_r^* Y_r^*$ of $\hat Y_r$.
		
		\textbf{Claim:}  There exists $R_1$ such that if $r \ge R_1$, then $\alpha_r \varphi_r \Gamma \subseteq \Gamma$.
		
		By applying Proposition \ref{contdis}.c to $\epsilon = p$, with probability one, there exists $R_2 > 0$ such that if $R \ge \max\{R_1, R_2\}$, then
			$$\nu_R( \eta_R \Theta_R g, B) \ge \nu_\oo( \alpha_R \eta_R \Theta_R g, B) - p \ge \inf_{y \in \Gamma} \nu_\oo(y, B) - p = p.$$

	\end{proof}

%% file: arxiv_dissertation.bbl
\begin{thebibliography}{LRST03}

\bibitem[AA09]{arguin2009structure}
L.P. Arguin and M.~Aizenman.
\newblock {On the structure of quasi-stationary competing particle systems}.
\newblock {\em Ann. Probab}, 37(3):1080--1113, 2009.

\bibitem[AL88]{arnold1988geometrical}
V.I. Arnolʹd and M.~Levi.
\newblock {\em {Geometrical methods in the theory of ordinary differential
  equations}}.
\newblock Springer, 1988.

\bibitem[Ale93]{alexander1993note}
K.S. Alexander.
\newblock {A note on some rates of convergence in first-passage percolation}.
\newblock {\em The Annals of Applied Probability}, 3(1):81--90, 1993.

\bibitem[Ale97]{alexander1997approximation}
K.S. Alexander.
\newblock {Approximation of subadditive functions and convergence rates in
  limiting-shape results}.
\newblock {\em The Annals of Probability}, 25(1):30--55, 1997.

\bibitem[ALR87]{aizenman1987some}
M.~Aizenman, J.L. Lebowitz, and D.~Ruelle.
\newblock {Some rigorous results on the Sherrington-Kirkpatrick spin glass
  model}.
\newblock {\em Communications in mathematical physics}, 112(1):3--20, 1987.

\bibitem[ASS07]{aizenman2007mean}
M.~Aizenman, R.~Sims, and S.L. Starr.
\newblock {Mean-Field Spin Glass models from the Cavity--ROSt Perspective}.
\newblock In {\em Prospects in mathematical physics: Young Researchers
  Symposium of the 14th International Congress on Mathematical Physics, July
  25-26, 2003, Lisbon, Portugal}, page~1. Amer Mathematical Society, 2007.

\bibitem[BDJ99]{baik1999distribution}
J.~Baik, P.~Deift, and K.~Johansson.
\newblock {On the distribution of the length of the longest increasing
  subsequence of random permutations}.
\newblock {\em Journal of the American Mathematical Society}, 12(4):1119--1178,
  1999.

\bibitem[BKS03]{benjamini2003first}
I.~Benjamini, G.~Kalai, and O.~Schramm.
\newblock {First passage percolation has sublinear distance variance}.
\newblock {\em Annals of Probability}, pages 1970--1978, 2003.

\bibitem[Bog98]{bogachev1998gm}
V.I. Bogachev.
\newblock {\em {Gaussian measures}}.
\newblock American Mathematical Society, 1998.

\bibitem[Boi90]{boivin1990first}
D.~Boivin.
\newblock {First passage percolation: the stationary case}.
\newblock {\em Probability Theory and Related Fields}, 86(4):491--499, 1990.

\bibitem[Bol89]{bolthausen1989note}
E.~Bolthausen.
\newblock {A note on the diffusion of directed polymers in a random
  environment}.
\newblock {\em Communications in Mathematical Physics}, 123(4):529--534, 1989.

\bibitem[BR06]{benaim2006modified}
M.~Bena\"im and R.~Rossignol.
\newblock {A modified Poincare inequality and its application to First Passage
  Percolation}.
\newblock {\em arXiv preprint math/0602496}, 2006.

\bibitem[BR08]{benaim2008exponential}
M.~Bena\"im and R.~Rossignol.
\newblock {Exponential concentration for First Passage Percolation through
  modified Poincar{\'e} inequalities}.
\newblock {\em Annales de l'Institut Henri Poincar{\'e}, Probabilit{\'e}s et
  Statistiques}, 44(3):544--573, 2008.

\bibitem[BS10]{blair2007first}
N.D. Blair-Stahn.
\newblock {First Passage Percolation and Competition Models}.
\newblock {\em arXiv preprint arXiv:1005.0649}, 2010.

\bibitem[CD81]{cox1981slt}
J.T. Cox and R.~Durrett.
\newblock Some limit theorems for percolation processes with necessary and
  sufficient conditions.
\newblock {\em The Annals of Probability}, 9(4):583--603, 1981.

\bibitem[CD09]{chatterjee2009central}
S.~Chatterjee and P.S. Dey.
\newblock {Central limit theorem for first-passage percolation time across thin
  cylinders}.
\newblock {\em arXiv preprint arXiv:0911.5702}, 2009.

\bibitem[CH02]{carmona2002partition}
P.~Carmona and Y.~Hu.
\newblock {On the partition function of a directed polymer in a Gaussian random
  environment}.
\newblock {\em Probability Theory and Related Fields}, 124(3):431--457, 2002.

\bibitem[Cha08]{chatterjee2008chaos}
S.~Chatterjee.
\newblock {Chaos, concentration, and multiple valleys}.
\newblock {\em arXiv preprint arXiv:0810.4221}, 2008.

\bibitem[Cha09]{chatterjee2009disorder}
S.~Chatterjee.
\newblock {Disorder chaos and multiple valleys in spin glasses}.
\newblock {\em Arxiv preprint arXiv:0907.3381}, 2009.

\bibitem[CSY03]{comets2003directed}
F.~Comets, T.~Shiga, and N.~Yoshida.
\newblock {Directed polymers in a random environment: path localization and
  strong disorder}.
\newblock {\em Bernoulli}, pages 705--723, 2003.

\bibitem[CSY04]{comets2004probabilistic}
F.~Comets, T.~Shiga, and N.~Yoshida.
\newblock {Probabilistic analysis of directed polymers in a random environment:
  a review}.
\newblock {\em Advanced Studies in Pure Mathematics}, 39:115--142, 2004.

\bibitem[Der85]{derrida1985generalization}
B.~Derrida.
\newblock {A generalization of the random energy model which includes
  correlations between energies}.
\newblock {\em Journal de Physique Lettres}, 46(9):401--407, 1985.

\bibitem[Der90]{derrida1990directed}
B.~Derrida.
\newblock {Directed polymers in a random medium}.
\newblock {\em Physica A: Statistical and Theoretical Physics}, 163(1):71--84,
  1990.

\bibitem[dH09]{den2009random}
F.~den Hollander.
\newblock {\em {Random polymers: {\'E}cole d'{\'E}t{\'e} de Probabilit{\'e}s de
  Saint-Flour XXXVII-2007}}.
\newblock Springer Verlag, 2009.

\bibitem[DL81]{durrett1981shape}
R.~Durrett and T.M. Liggett.
\newblock {The shape of the limit set in Richardson's growth model}.
\newblock {\em The Annals of Probability}, 9(2):186--193, 1981.

\bibitem[Dur96]{durrett1996probability}
R.~Durrett.
\newblock {\em {Probability: theory and examples}}.
\newblock Duxbury Press Belmont, CA, 1996.

\bibitem[EA75]{edwards1975theory}
S.F. Edwards and P.W. Anderson.
\newblock {Theory of spin glasses}.
\newblock {\em Journal of Physics F: Metal Physics}, 5:965--974, 1975.

\bibitem[Fol99]{folland1999real}
G.B. Folland.
\newblock {\em {Real Analysis: Modern Techniques and Their Applications}}.
\newblock Wiley-Interscience, 1999.

\bibitem[GH75]{geman1975random}
D.~Geman and J.~Horowitz.
\newblock {Random shifts which preserve measure}.
\newblock {\em Proceedings of the American Mathematical Society},
  49(1):143--150, 1975.

\bibitem[Gia07]{giacomin2007random}
G.~Giacomin.
\newblock {\em {Random polymer models}}.
\newblock Imperial College Pr, 2007.

\bibitem[GT02]{guerra2002thermodynamic}
F.~Guerra and F.L. Toninelli.
\newblock {The thermodynamic limit in mean field spin glass models}.
\newblock {\em Communications in Mathematical Physics}, 230(1):71--79, 2002.

\bibitem[GT06]{giacomin2006smoothing}
G.~Giacomin and F.L. Toninelli.
\newblock {Smoothing effect of quenched disorder on polymer depinning
  transitions}.
\newblock {\em Communications in Mathematical Physics}, 266(1):1--16, 2006.

\bibitem[HH85]{huse1985pinning}
D.A. Huse and C.L. Henley.
\newblock {Pinning and roughening of domain walls in Ising systems due to
  random impurities}.
\newblock {\em Physical review letters}, 54(25):2708--2711, 1985.

\bibitem[HHF85]{huse1985huse}
D.A. Huse, C.L. Henley, and D.S. Fisher.
\newblock {Huse, Henley, and Fisher respond}.
\newblock {\em Physical Review Letters}, 55(26):2924--2924, 1985.

\bibitem[HM95]{haggstrom1995ass}
O.~H{\"a}ggstr{\"o}m and R.~Meester.
\newblock Asymptotic shapes for stationary first passage percolation.
\newblock {\em The Annals of Probability}, 23(4):1511--1522, 1995.

\bibitem[HM07]{hambly2007heavy}
B.~Hambly and J.B. Martin.
\newblock {Heavy tails in last-passage percolation}.
\newblock {\em Probability Theory and Related Fields}, 137(1):227--275, 2007.

\bibitem[HN97]{howard1997euclidean}
C.D. Howard and C.M. Newman.
\newblock {Euclidean models of first-passage percolation}.
\newblock {\em Probability Theory and Related Fields}, 108(2):153--170, 1997.

\bibitem[HN00]{howard2000geodesics}
C.D. Howard and C.M. Newman.
\newblock {Geodesics and spanning trees for Euclidean first-passage
  percolation}.
\newblock {\em arXiv preprint math/0010205}, 2000.

\bibitem[HN01]{howard2001special}
C.D. Howard and C.M. Newman.
\newblock {Special Invited Paper: Geodesics And Spanning Trees For Euclidean
  First Passage Percolation}.
\newblock {\em Ann. Probab}, 29(2):577--623, 2001.

\bibitem[Hof08]{hoffman2008geodesics}
C.~Hoffman.
\newblock {Geodesics in first passage percolation}.
\newblock {\em Ann. Appl. Probab}, 18(5):1944--1969, 2008.

\bibitem[How00]{howard2000lower}
C.D. Howard.
\newblock {Lower bounds for point-to-point wandering exponents in Euclidean
  first-passage percolation}.
\newblock {\em Journal of Applied Probability}, 37(4):1061--1073, 2000.

\bibitem[How04]{howard2004mfp}
C.D. Howard.
\newblock Models of first-passage percolation.
\newblock {\em Probability on Discrete Structures}, pages 125--173, 2004.

\bibitem[HP98]{häggström1998first}
O.~H{\"a}ggstr{\"o}m and R.~Pemantle.
\newblock {First passage percolation and a model for competing spatial growth}.
\newblock {\em Journal of Applied Probability}, 35(3):683--692, 1998.

\bibitem[HW65]{hammersley1965fpp}
J.M. Hammersley and D.J.A. Welsh.
\newblock {First-passage percolation, sub-additive process, stochastic network
  and generalized renewal theory}.
\newblock {\em Bernoulli, 1713: Bayes, 1763; Laplace, 1813. Anniversary
  Volume}, page~61, 1965.

\bibitem[IR78]{ibragimov1978grp}
I.A. Ibragimov and J.A. Rozanov.
\newblock {\em Gaussian random processes}.
\newblock Springer-Verlag New York, 1978.

\bibitem[IS88]{imbrie1988diffusion}
J.Z. Imbrie and T.~Spencer.
\newblock {Diffusion of directed polymers in a random environment}.
\newblock {\em Journal of Statistical Physics}, 52(3):609--626, 1988.

\bibitem[Joh00]{johansson2000shape}
K.~Johansson.
\newblock {Shape fluctuations and random matrices}.
\newblock {\em Communications in mathematical physics}, 209(2):437--476, 2000.

\bibitem[Kes84]{kesten1180arp}
H.~Kesten.
\newblock Aspects of first passage percolation.
\newblock {\em Ecole d'\'et\'e de Probabilit\'es de St. Flour. Lecture Notes in
  Math}, 1180:125--264, 1984.

\bibitem[Kes93]{kesten1993scf}
H.~Kesten.
\newblock On the speed of convergence in first-passage percolation.
\newblock {\em The Annals of Applied Probability}, 3(2):296--338, 1993.

\bibitem[Kin68]{kingman1968ets}
J.F.C. Kingman.
\newblock The ergodic theory of subadditive stochastic processes.
\newblock {\em Journal of the Royal Statistical Society. Series B
  (Methodological)}, 30(3):499--510, 1968.

\bibitem[KN85]{kardar1985commensurate}
M.~Kardar and D.R. Nelson.
\newblock {Commensurate-incommensurate transitions with quenched random
  impurities}.
\newblock {\em Physical review letters}, 55(11):1157--1160, 1985.

\bibitem[KPZ86]{kardar1986dynamic}
M.~Kardar, G.~Parisi, and Y.C. Zhang.
\newblock {Dynamic scaling of growing interfaces}.
\newblock {\em Physical Review Letters}, 56(9):889--892, 1986.

\bibitem[KS88]{krug1988universality}
J.~Krug and H.~Spohn.
\newblock {Universality classes for deterministic surface growth}.
\newblock {\em Physical Review A}, 38(8):4271--4283, 1988.

\bibitem[KS91]{krug1991kinetic}
J.~Krug and H.~Spohn.
\newblock {Kinetic roughening of growing surfaces}.
\newblock {\em Solids far from equilibrium}, pages 479--582, 1991.

\bibitem[LaG10]{lagatta2010continuous}
T.~LaGatta.
\newblock {Continuous Disintegrations of Gaussian Measures}.
\newblock {\em arXiv preprint arXiv:1003.0975}, 2010.

\bibitem[Lee97]{lee1997rmi}
J.M. Lee.
\newblock {\em Riemannian Manifolds: An Introduction to Curvature}.
\newblock Springer, 1997.

\bibitem[LN96]{licea1996gtd}
C.~Licea and C.M. Newman.
\newblock Geodesics in two-dimensional first-passage percolation.
\newblock {\em The Annals of Probability}, 24(1):399--410, 1996.

\bibitem[LNP96]{licea1996superdiffusivity}
C.~Licea, C.M. Newman, and M.S.T. Piza.
\newblock {Superdiffusivity in first-passage percolation}.
\newblock {\em Probability Theory and Related Fields}, 106(4):559--591, 1996.

\bibitem[LRST03]{lamburt2003grc}
V.G. Lamburt, E.R. Rozendorn, D.D. Sokoloff, and V.N. Tutubalin.
\newblock {Geodesics with random curvature on Riemannian and pseudo-Riemannian
  manifolds}.
\newblock {\em Trudy Geom. Sem. Kazan Gos. Univ}, 24:99--106, 2003.

\bibitem[LW10]{lagatta2009shape}
T.~LaGatta and J.~Wehr.
\newblock {A Shape Theorem for Riemannian First-Passage Percolation}.
\newblock {\em J. Math. Phys.}, 51(5), 2010.

\bibitem[Mar00]{martin2000linear}
J.B. Martin.
\newblock {Linear growth for greedy lattice animals}.
\newblock {\em Rapport de recherche-institut national de recherche en
  informatique et en automatique}, 2000.

\bibitem[Mar04]{martin2004limiting}
J.B. Martin.
\newblock {Limiting shape for directed percolation models}.
\newblock {\em Annals of probability}, 32(4):2908--2937, 2004.

\bibitem[MPV87]{mezard1987spin}
M.~Mezard, G.~Parisi, and M.A. Virasoro.
\newblock {\em {Spin glass theory and beyond}}.
\newblock World Scientific Singapore, 1987.

\bibitem[New95]{newman1995svf}
C.M. Newman.
\newblock A surface view of first-passage percolation.
\newblock In S.D. Chatterji, editor, {\em Proceedings of the International
  Congress of Mathematicians}, volume~2, pages 1017--1023, 1995.

\bibitem[New97]{newman1997tds}
C.M. Newman.
\newblock {\em Topics in disordered systems}.
\newblock Birkh{\"a}user, 1997.

\bibitem[New10]{newmanemailapr22}
C.M. Newman.
\newblock personal communication, 2010.

\bibitem[NP95]{newman1995dsf}
C.M. Newman and M.S.T. Piza.
\newblock Divergence of shape fluctuations in two dimensions.
\newblock {\em The Annals of Probability}, 23(3):977--1005, 1995.

\bibitem[Par79]{parisi1979infinite}
G.~Parisi.
\newblock {Infinite number of order parameters for spin-glasses}.
\newblock {\em Physical Review Letters}, 43(23):1754--1756, 1979.

\bibitem[Piz97]{piza1997directed}
M.S.T. Piza.
\newblock {Directed polymers in a random environment: some results on
  fluctuations}.
\newblock {\em Journal of Statistical Physics}, 89(3):581--603, 1997.

\bibitem[PRT00]{parisi2000origin}
G.~Parisi and F.~Ricci-Tersenghi.
\newblock {On the origin of ultrametricity}.
\newblock {\em Journal of Physics A: Mathematical and General}, 33:113--129,
  2000.

\bibitem[Ric73]{richardson1973random}
D.~Richardson.
\newblock {Random growth in a tessellation}.
\newblock In {\em Proceedings of the Cambridge Philosophical Society},
  volume~74, page 515, 1973.

\bibitem[Rue87]{ruelle1987mathematical}
D.~Ruelle.
\newblock {A mathematical reformulation of Derrida's REM and GREM}.
\newblock {\em Communications in Mathematical Physics}, 108(2):225--239, 1987.

\bibitem[SK75]{sherrington1975solvable}
D.~Sherrington and S.~Kirkpatrick.
\newblock {Solvable model of a spin-glass}.
\newblock {\em Physical review letters}, 35(26):1792--1796, 1975.

\bibitem[SZ96]{song1996remark}
R.~Song and X.Y. Zhou.
\newblock {A remark on diffusion of directed polymers in random environments}.
\newblock {\em Journal of Statistical Physics}, 85(1):277--289, 1996.

\bibitem[Tal94]{talagrand1994russo}
M.~Talagrand.
\newblock {On Russo's approximate zero-one law}.
\newblock {\em The Annals of Probability}, 22(3):1576--1587, 1994.

\bibitem[Tal98]{talagrand1998sherrington}
M.~Talagrand.
\newblock {The Sherrington--Kirkpatrick model: a challenge for mathematicians}.
\newblock {\em Probability Theory and Related Fields}, 110(2):109--176, 1998.

\bibitem[TW94]{tracy1994level}
C.A. Tracy and H.~Widom.
\newblock {Level-spacing distributions and the Airy kernel}.
\newblock {\em Communications in Mathematical Physics}, 159(1):151--174, 1994.

\bibitem[VAW90]{vahidi1990first}
M.Q. Vahidi-Asl and J.C. Wierman.
\newblock {First-passage percolation on the Voronoi tessellation and Delaunay
  triangulation}.
\newblock In {\em Random graphs' 87: based on proceedings of the 3rd
  International Seminar on Random Graphs and Probabilistic Methods in
  Combinatorics, June 27-July 3 1987, Pozna\'n, Poland}, page 341. John Wiley
  \& Sons Inc, 1990.

\bibitem[VAW92]{vahidi1992shape}
M.Q. Vahidi-Asl and J.C. Wierman.
\newblock {A shape result for first-passage percolation on the Voronoi
  tessellation and Delaunay triangulation}.
\newblock In {\em Random graphs}, volume~2, pages 247--262. Wiley-Interscience,
  1992.

\bibitem[WA90]{wehr1990fluctuations}
J.~Wehr and M.~Aizenman.
\newblock {Fluctuations of extensive functions of quenched random couplings}.
\newblock {\em Journal of Statistical Physics}, 60(3):287--306, 1990.

\bibitem[Weh97]{wehr1997number}
J.~Wehr.
\newblock {On the number of infinite geodesics and ground states in disordered
  systems}.
\newblock {\em Journal of Statistical Physics}, 87(1):439--447, 1997.

\bibitem[WW98]{wehr1998absence}
J.~Wehr and J.~Woo.
\newblock {Absence of geodesics in first-passage percolation on a half-plane}.
\newblock {\em Annals of Probability}, 26(1):358--367, 1998.

\bibitem[Zir01]{zirbel2001lagrangian}
C.L. Zirbel.
\newblock {Lagrangian observations of homogeneous random environments}.
\newblock {\em Advances in Applied Probability}, 33(4):810--835, 2001.

\end{thebibliography}
